\documentclass[12pt]{article}
\usepackage{amssymb}
\usepackage{amsfonts}
\usepackage{latexsym}
\usepackage{amsthm}

\newcommand{\er}[1]{{\rm(\ref{#1})}}
\def\lb{\label}
\theoremstyle{plain}
\newtheorem{theorem}{\bf Theorem}[section]
\newtheorem{lemma}[theorem]{\bf Lemma}

\newtheorem{proposition}[theorem]{\bf Proposition}
\theoremstyle{remark}

\begin{document}

\def\a{\alpha} \def\cA{{\cal A}} \def\bA{{\bf A}}  \def\mA{{\mathscr A}}
\def\b{\beta}  \def\cB{{\cal B}} \def\bB{{\bf B}}  \def\mB{{\mathscr B}}
\def\g{\gamma} \def\cC{{\cal C}} \def\bC{{\bf C}}  \def\mC{{\mathscr C}}
\def\G{\Gamma} \def\cD{{\cal D}} \def\bD{{\bf D}}  \def\mD{{\mathscr D}}
\def\d{\delta} \def\cE{{\cal E}} \def\bE{{\bf E}}  \def\mE{{\mathscr E}}
\def\D{\Delta} \def\cF{{\cal F}} \def\bF{{\bf F}}  \def\mF{{\mathscr F}}
\def\c{\chi}   \def\cG{{\cal G}} \def\bG{{\bf G}}  \def\mG{{\mathscr G}}
\def\z{\zeta}  \def\cH{{\cal H}} \def\bH{{\bf H}}  \def\mH{{\mathscr H}}
\def\e{\eta}   \def\cI{{\cal I}} \def\bI{{\bf I}}  \def\mI{{\mathscr I}}
\def\p{\psi}   \def\cJ{{\cal J}} \def\bJ{{\bf J}}  \def\mJ{{\mathscr J}}
\def\vT{\Theta}\def\cK{{\cal K}} \def\bK{{\bf K}}  \def\mK{{\mathscr K}}
\def\k{\kappa} \def\cL{{\cal L}} \def\bL{{\bf L}}  \def\mL{{\mathscr L}}
\def\l{\lambda}\def\cM{{\cal M}} \def\bM{{\bf M}}  \def\mM{{\mathscr M}}
\def\L{\Lambda}\def\cN{{\cal N}} \def\bN{{\bf N}}  \def\mN{{\mathscr N}}
\def\m{\mu}    \def\cO{{\cal O}} \def\bO{{\bf O}}  \def\mO{{\mathscr O}}
\def\n{\nu}    \def\cP{{\cal P}} \def\bP{{\bf P}}  \def\mP{{\mathscr P}}
\def\r{\rho}   \def\cQ{{\cal Q}} \def\bQ{{\bf Q}}  \def\mQ{{\mathscr Q}}
\def\s{\sigma} \def\cR{{\cal R}} \def\bR{{\bf R}}  \def\mR{{\mathscr R}}
\def\S{\Sigma} \def\cS{{\cal S}} \def\bS{{\bf S}}  \def\mS{{\mathscr S}}
\def\t{\tau}   \def\cT{{\cal T}} \def\bT{{\bf T}}  \def\mT{{\mathscr T}}
\def\f{\phi}   \def\cU{{\cal U}} \def\bU{{\bf U}}  \def\mU{{\mathscr U}}
\def\F{\Phi}   \def\cV{{\cal V}} \def\bV{{\bf V}}  \def\mV{{\mathscr V}}
\def\P{\Psi}   \def\cW{{\cal W}} \def\bW{{\bf W}}  \def\mW{{\mathscr W}}
\def\o{\omega} \def\cX{{\cal X}} \def\bX{{\bf X}}  \def\mX{{\mathscr X}}
\def\x{\xi}    \def\cY{{\cal Y}} \def\bY{{\bf Y}}  \def\mY{{\mathscr Y}}
\def\X{\Xi}    \def\cZ{{\cal Z}} \def\bZ{{\bf Z}}  \def\mZ{{\mathscr Z}}
\def\O{\Omega}
\def\ve{\varepsilon}
\def\vt{\vartheta}
\def\vp{\varphi}
\def\vk{\varkappa}

\def\mM{\cM}
\def\mB{\cB}
\def\mR{\cR}

\def\Z{{\Bbb Z}}
\def\R{{\Bbb R}}
\def\C{{\Bbb C}}
\def\T{{\Bbb T}}
\def\N{{\Bbb N}}
\def\S{{\Bbb S}}
\def\H{{\Bbb H}}

\let\ge\geqslant
\let\le\leqslant
\let\geq\geqslant
\let\leq\leqslant
\def\ma{\left(\begin{array}{cc}}
\def\am{\end{array}\right)}
\def\iint{\int\!\!\!\int}
\def\lt{\biggl}
\def\rt{\biggr}
\let\geq\geqslant
\let\leq\leqslant
\def\[{\begin{equation}}
\def\]{\end{equation}}
\def\wt{\widetilde}
\def\pa{\partial}
\def\sm{\setminus}
\def\es{\emptyset}
\def\no{\noindent}
\def\ol{\overline}
\def\iy{\infty}
\def\ev{\equiv}
\def\/{\over}
\def\ts{\times}
\def\os{\oplus}
\def\ss{\subset}
\def\h{\hat}
\def\Re{\mathop{\rm Re}\nolimits}
\def\Im{\mathop{\rm Im}\nolimits}
\def\supp{\mathop{\rm supp}\nolimits}
\def\sign{\mathop{\rm sign}\nolimits}
\def\Ran{\mathop{\rm Ran}\nolimits}
\def\Ker{\mathop{\rm Ker}\nolimits}
\def\Tr{\mathop{\rm Tr}\nolimits}
\def\const{\mathop{\rm const}\nolimits}
\def\dist{\mathop{\rm dist}\nolimits}
\def\diag{\mathop{\rm diag}\nolimits}
\def\Wr{\mathop{\rm Wr}\nolimits}
\def\BBox{\hspace{1mm}\vrule height6pt width5.5pt depth0pt \hspace{6pt}}

\def\Diag{\mathop{\rm Diag}\nolimits}

\def\Twelve{
\font\Tenmsa=msam10 scaled 1200
\font\Sevenmsa=msam7 scaled 1200
\font\Fivemsa=msam5 scaled 1200
\textfont\msbfam=\Tenmsb
\scriptfont\msbfam=\Sevenmsb
\scriptscriptfont\msbfam=\Fivemsb

\font\Teneufm=eufm10 scaled 1200
\font\Seveneufm=eufm7 scaled 1200
\font\Fiveeufm=eufm5 scaled 1200
\textfont\eufmfam=\Teneufm
\scriptfont\eufmfam=\Seveneufm
\scriptscriptfont\eufmfam=\Fiveeufm}

\def\Ten{
\textfont\msafam=\tenmsa
\scriptfont\msafam=\sevenmsa
\scriptscriptfont\msafam=\fivemsa

\textfont\msbfam=\tenmsb
\scriptfont\msbfam=\sevenmsb
\scriptscriptfont\msbfam=\fivemsb

\textfont\eufmfam=\teneufm
\scriptfont\eufmfam=\seveneufm
\scriptscriptfont\eufmfam=\fiveeufm}

\title {Spectral asymptotics for periodic fourth order operators}

\author{Andrei Badanin
\begin{footnote}
{ Department of  Mathematics of Archangel University, Russia e-mail: a.badanin@agtu.ru}
\end{footnote}
 \and Evgeny Korotyaev
\begin{footnote}
{ Institut f\"ur  Mathematik,  Humboldt Universit\"at zu Berlin,
Rudower Chaussee 25, 12489, Berlin, Germany,
e-mail: evgeny@math.hu-berlin.de, 
corresponding author}
\end{footnote}
}

\maketitle

\begin{abstract}
\no We consider the operator ${d^4\/dt^4}+V$ on the real line
with a real periodic potential $V$. 
The spectrum of this operator is absolutely continuous and consists of intervals  separated by gaps. We define a Lyapunov function
which is analytic on a two sheeted Riemann surface. On each
sheet, the Lyapunov function has the same properties 
as in the scalar case, but it has 
branch points, which we call resonances. We prove the existence of real as well as non-real resonances for specific potentials. We determine the asymptotics of the periodic and anti-periodic spectrum and of the resonances at high energy. We show that there exist two type of
gaps: 1) stable gaps, where the endpoints are periodic and
anti-periodic eigenvalues, 2) unstable (resonance) gaps, where the
endpoints are resonances (i.e., real branch points of the Lyapunov
function above the bottom of the spectrum). We also show that
the periodic and anti-periodic spectrum together determine the spectrum of our operator. Finally, we show that for small potentials $V\ne 0$ the spectrum in the lowest band  has multiplicity 4 and the bottom of the spectrum is a resonance, and not a periodic (or anti-periodic) eigenvalue.
\end{abstract}


\section {Introduction and main results}
\setcounter{equation}{0}

We consider the self-adjoint operator $\cL={d^4\/dt^4}+V,$ acting on
$L^2(\R)$, where the real 1-periodic potential $V$ belongs to the
real space $L^1_0(\T)=\{V\in L^1(\T),\int_0^1V(t)dt=0\},\  \T=\R/\Z$, equipped with the norm $\|V\|=\int_0^1|V(t)|dt<\iy$.
It is well known (see [DS]) that the spectrum $\s(\cL)$ of $\cL$ is
absolutely continuous and consists of non-degenerate intervals.
These intervals are separated by the gaps $G_n=(E_n^-,E_n^+), n\ge 1$, with length $|G_n|>0$. Introduce the fundamental solutions
$\vp_j(t,\l), j=0,1,2,3,$  of the equation
\[
\lb{1b}
y''''+Vy=\l y,\ \ \ \ \ \ (t,\l)\in \R \ts \C,
\]
satisfying the following conditions: $\vp_j^{(k)}(0,\l)=\d_{jk},j,k=0,...,3$, where $\d_{jk}$
is the standard Kronecker symbol.
Here and below we use the notation $f'={\pa f\/\pa t},
f^{(k)}={\pa^kf \/\pa t^k}$.
We define the monodromy $4\ts 4$-matrix $M$ by
\[
\lb{1c} M(\l)=\mM(1,\l),\ \ \ \mM(t,\l)
=\{\mM_{kj}(t,\l)\}_{j,k=0}^3
=\{\vp_j^{(k)}(t,\l)\}\}_{j,k=0}^3.
\]
The matrix valued function $M$ is entire. An eigenvalue of $M(\l)$
is called a {\it multiplier}. It is a root of the algebraic
equation $D(\t,\l)=0$, where $D(\t,\l)\ev\det(M(\l)-\t I_4),
\t,\l\in\C$. Let $D_{\pm}(\l)={1\/4}D(\pm 1,\l)$. The zeros of
$D_{+}(\l)$ (or $D_{-}(\l)$) are the eigenvalues of the periodic
(anti-periodic) problem for the equation $y''''+Vy=\l y$. Denote by
$\l_0^+,\l_{2n}^\pm, n=1,2,...$ the sequence of zeros of $D_+$
(counted with multiplicity) such that $\l_{0}^+\le \l_{2}^-\le
\l_{2}^+\le \l_{4}^-\le\l_{4}^+\le \l_{6}^- \le... $ Denote by
$\l_{2n-1}^\pm, n=1,2,...$ the sequence of zeros of $D_-$ (counted
with multiplicity) such that $\l_{1}^-\le \l_{1}^+\le \l_{3}^-\le
\l_{3}^+\le\l_{5}^-\le \l_{5}^+ \le... $.

A great number of  papers is devoted to the inverse spectral
theory for the Hill operator. We mention all papers where the
inverse problem including characterization was solved: Marchenko
and Ostrovski [MO], Garnett and Trubowitz [GT1-2],Kappeler [Kap],
Kargaev and Korotyaev [KK1], and Korotyaev [K1-3] and for $2\ts 2$
Dirac operator [Mi1-2], [K4-5]. Recently, one of the authors [K6]
extended the results of [MO], [GT1], [K1-2] of the case $-y''+uy
$ to the case of distributions, i.e. $-y''+u'y$ on $L^2(\R)$,
where $u\in L_{loc}^2(\R)$ is periodic.

There exist many papers about the periodic systems $N\ge 2$ (see [YS]). The basic results for the direct spectral theory for the matrix case were obtained by Lyapunov [Ly]  (see also the interesting papers of Krein [Kr], Gel'fand and Lidskii [GL]).
 The operator $-{d^2\/dt^2}+\cV$ on the real line where $\cV$  is a 1-periodic $2\ts 2$ matrix potential was considered in [BBK]. The following results are obtained: the Lyapunov function is constructed as an analytic function on a 2-sheeted 
Riemann surface and the existence of real and complex resonances
are proved for some specific potentials. 
Recall the well-known Lyapunov Theorem,
 in a formulation adapted for our case  (see [Ly],[YS]):

{\bf Theorem (Lyapunov)} {\it Let $V\in L_0^1(\T)$. Then
\[
\lb{rD} D(\t,\cdot)=\t^4 D(\t^{-1},\cdot),\ \ \ \t\ne 0.
\]
If for some $\l\in\C$ (or $\l\in\R$) $\t(\l)$  is a multiplier of multiplicity $d\ge 1$, then $\t^{-1}(\l)$ (or $\ol\t(\l)$) is a multiplier of multiplicity $d$.
Moreover, each $M(\l), \l\in\C$ has exactly four multipliers
$\t_1^{\pm 1}(\l), \t_2^{\pm 1}(\l)$. Furthermore, $\l\in\s(\cL)$ iff $|\t_1(\l)|=1 $ or $|\t_2(\l)|=1$.
 If $\t(\l)$ is a simple multiplier and $|\t(\l)|=1$, then
$\t'(\l)\ne 0$.}

The spectral problems for the fourth order periodic operator were the subject of many authors (see [DS],[GO],[MV],[P1-2],[PK1-2],[YS]). Firstly, we mention the papers of Papanicolaou [P1-2] devoted to the  Euler-Bernoulli
equation $(ay'')''=\l by$ with the periodic functions $a,b$.
For this case he defines the Lyapunov function and obtains some properties of this function. In particular, it is proved that the Lyapunov function is analytic on some two sheeted Riemann surface. It is important that for this case he proved that all branch points of the Lyapunov function are real and $\le 0$.
Note that in our case we have the example, Proposition \ref{1.4}, where the Lyapunov function has real and  non-real branch points. This the main difference between our Lyapunov function
 and his one.
Moreover, Papanicolaou proved that if all the gaps are closed and the Lyapunov function is entire in $\sqrt \l$, then the functions $a,b$ are constants.
Secondly, Papanicolaou and Kravvaritis [PK1-2] considered some  inverse problems for the  Euler-Bernoulli equation.
Thirdly,  Galunov, Oleinik [GO] considered the operator $y^{(2n)}+\g \d_{per}(t)y,\g\in\R$, on the real line  with the
periodic delta-potential $\d_{per}(t)=\sum_{n\in \Z} \d(t-n)$. They study the spectrum in the lowest band. It can be seen from the results of this paper, that the spectrum in this band
has multiplicity $4$, for $n=2$ and for some $\g$.
Recall that  (see Theorem \ref{1.3}) we prove a stronger result and show that the spectrum has multiplicity $4$ in the lowest band for each sufficiently small potential.

We introduce the functions
\[
\lb{1d} T_m={1\/4}\Tr M^m,\ \ \ m\ge 1,\ \ \
\r={T_2+1\/2}-T_1^2.
\]
The functions $T_1,T_2,\r$ are real on $\R$ and entire.
In the case  $V=0$ the corresponding functions have the forms
\[
\lb{0} T_m^0={\cos mz+\cosh mz\/2},\ \ m=1,2, \ \r^0={(\cos z-\cosh
z)^2\/4}, \ \ \ z=\l^{1\/4},
\]
here and below $\arg z\in (-{\pi\/4},{\pi\/4}]$. It is known (see [RS])  that $D(\t,\l)=\sum_0^4\x_m(\l)\t^{4-m}$, where the functions $\x_m$ are given by
$$
\x_0=1,\ \ \ \x_1=-4T_1,\ \ \ \x_2=-2(T_2+T_1\x_1),\ \  . \ . \ .\ \ .
$$
 Then using the identity \er{rD} we obtain $D(\t,\cdot)=(\t^4+1)+\x_1(\t^3+\t)+\x_2\t^2$,
which yields
\[
\lb{2l} D(\t,\cdot)
=\lt(\t^2-2(T_1-\sqrt{\r})\t+1\rt)
\lt(\t^2-2(T_1+\sqrt{\r})\t+1\rt).
\]
We introduce the domains
\[
\lb{DD}
\cD_r=\rt\{\l\in \C: |\l^{1/4}|>r,\ \  |\l^{1/4}-(1\pm i)\pi n|>{\pi\/4},|\l^{1/4}-\pi n|>{\pi\/4},n\ge 0 \rt\},\ \ \ r\ge 0.
\]
We have $\r(\l)=\r^0(\l)(1+o(1))$
as $|\l|\to \iy,\l\in\cD_1$ (see Lemma \ref{ar}). Then we define the analytic
function $\sqrt{\r(\l)},\l\in\cD_r
$ for some large $r>0$, by the condition
$\sqrt{\r(\l)}=\sqrt{\r^0(\l)}(1+o(1))$ as $|\l|\to\iy,\l\in \cD_r$,
where $\sqrt{\r^0(\l)}={(\cos z-\cosh
z)/2}$.

The function $\r$ is real on $\R$, then $r$ is a root of $\r$ iff $\ol r$ is a root of $\r$.  
By Lemma 5.1, for large integer $N$ the function $\r(\l)$ has exactly $2N+1$ roots, counted with multiplicity, in the disk
$\{\l:|\l|<4(\pi(N+{1\/2}))^4\}$ and for each $n>N$, exactly two
roots, counted with multiplicity, in the domain
$\{\l:|\l^{1/4}-\pi(1+i)n|<\pi/4\}$. There are no other roots.
Thus the function $\r(\l)$ has an odd number $\ge 1$ of real zeros
(counted with multiplicity) on the real interval 
$(-\G, \G)\ss \R, \G=4(\pi(N+{1\/2}))^4$.

Let $\{r_0^-,r_{n}^\pm\}_{1}^{\iy }$ be the sequence of
 zeros of $\r$ in $\C$ (counted with multiplicity) such that:

$r_0^-$ is the maximal real zero, and $..\le \Re r_{n+1}^+\le \Re r_{n}^+\le ... \le \Re r_{1}^+$,

if $r_{n}^+\in \C_+$, then $r_{n}^-=\ol{r_{n}^+}\in \C_-$,

if $r_{n}^+\in \R$, then $r_{n}^-\le r_{n}^+\le \Re r_{m-1}^-, m=1,..,n$.

Below we will show that $r_n^\pm =-4(\pi n)^4+0(n^2)$ as $n\to \iy$,
see Lemma 5.1.
We call a zero of $\r$ a resonance of $\cL$.
Let $...\le r_{n_j}^-\le r_{n_{j}}^+\le...\le r_{n_1}^-\le r_{n_1}^+\le r_{0}^-$
be the subsequence of the real zeros of $\r$.
Then  $\r(\l)<0$  for any $\l\in \g_j^0=(r_{n_{j+1}}^+,r_{n_j}^-),j\ge 1$.
We call an interval $\g_j^0\ss \R$ a resonance gap.

We construct the Riemann surface $\mR$ for $\sqrt \r$.
 For any $r_n^+\in\C_+$ we take some curve $\e_n$, which joins the
points $r_n^+, \ol{r_n^+}$ and does not cross $\g^0=\cup \g_n^0$.
To "build" the surface $\mR$, we take two
replicas of the $\l$-plane cut along $\g^0$ and $\cup \e_n$ and call
them sheet $\mR_1$ and sheet $\mR_2$. The cut on each sheet
has two edges; we label each edge with a + or a $-$. Then
attach the $-$ edge of the cut on $\mR_1$ to the $+$ edge of
the cut on $\mR_2$ and attach the $+$ edge of the cut on
$\mR_1$ to the $-$ edge of the cut on $\mR_2$. Thus,
whenever we cross the cut, we pass from one sheet to the other.
 There exists a unique analytic continuation of the function
 $\sqrt \r$ from $\cD_r$ into the two sheeted Riemann surface
$\mR$ of the function $\sqrt \r$.
 Let below  $\z\in \mR$ and let $\f(\z)=\l$ be the natural projection $\f:\mR\to\C$.

 We introduce the Lyapunov function by
\[
\lb{1g}
 \D(\z)=T_1(\z)+\sqrt{\r(\z)},\ \ \ \ \ \ \ \z\in\mR.
\]
Note that $T_1(\z)=T_1(\l)$, since $T_1$ is entire.
Let $\D(\z)=\D_1(\l)$ on the first sheet $\cR_1$, $\D(\z)=\D_2(\l)$ on the second sheet
$\cR_2$. Then
\[
\lb{drs2} \D_1(\l)=T_1(\l)+\sqrt{\r(\l)},\ \
\D_2(\l)=T_1(\l)-\sqrt{\r(\l)},\ \ \l=\f(\z).
\]
Now we formulate our first result about the function
$\D(\l)$.

\begin{theorem}   \lb{T1}
 Let $V\in L_0^1(\T)$. Then the function
$\D=T_1+\sqrt{\r}$ is analytic on the two sheeted Riemann
surface $\mR$ and the branches $\D_m$ of
$\D$ have  the forms
\[  \lb{1ga}\lb{tD}
\D_m(\l)={\t_m(\l)+\t_m^{-1}(\l)\/2},
\ \ \ \l\in\mR_m,\ \ \ m=1,2,
\]
and the following properties:

\no i) The following identities and asymptotics are fulfilled
\[
\lb{1ca}\lb{aev}
\t_1(\l)=e^{\l^{1/4}+O(\l^{-{3/2}})},\ \ \
\t_2(\l)=e^{i\l^{1/4}+O(\l^{-{3/2}})},
\]
\[
\lb{asL}\lb{aD} \D_1(\l)=\cosh \l^{1/4}\lt(1+O(\l^{-3/2})\rt),\ \ \
\D_2(\l)=\cos\l^{1/4}\lt(1+O(\l^{-3/2})\rt)
\]
as $ |\l|\to\iy, \l\in\cD_1$.

\no ii) $\l\in \s(\cL)$ iff $\D_m(\l)\in [-1,1]$
for some $m=1,2$. Moreover, if $\l\in\s(\cL)$, then $\r(\l)\ge 0$.

\no iii) The spectrum of $\cL$ on an interval $S\ss \R$ has multiplicity
4 iff $-1<\D_m(z)<1$ for all $m=1,2, \l\in S$, except for finite
number of points.

\no  iv) the spectrum of $\cL$ on an interval $S\ss \R$ has multiplicity
2 iff $-1<\D_1(\l)<1, \D_2(\l)\in \R\sm [-1,1]$ or
$-1<\D_2(\l)<1, \D_1(\l)\in \R\sm [-1,1]$ 
for all $\l\in S$,except for finite number of points.

\no  v) Let $\D_m$ be real analytic on some interval
$I=(\a_1,\a_2)\ss\R$ and $-1<\D_m(\l)<1$, for any $\l\in I$ for
some $m\in \{1,2\}$. Then $\D_m'(\l)\ne 0$ for each $\l\in I$ (the
monotonicity property).

\no  vi) Each gap $G_n=(E_n^-,E_n^+), n\ge 1$ is a bounded interval and $E_n^\pm$ are either periodic (anti-periodic) eigenvalues or real branch point of $\D_m$ (for some $m=1,2$)
which is a zero of $\r$ (that is a resonance).
\end{theorem}

\no {\bf Remark.}
i) In the case of the Hill operator the monodromy matrix
has exactly 2 eigenvalues $\t,\t^{-1}$. The Lyapunov function
${1\/2}(\t+\t^{-1})$  is an entire function of the spectral
parameter. It defines the band structure of the spectrum. By
Theorem 1.1, the Lyapunov function for the operator $\cL$ also defines the band structure of the spectrum, but
it is an analytic function on a 2-sheeted Riemann surface.
The qualitative behavior of the Lyapunov function for small potentials is shown on Fig. \ref{sp}.

\no ii) Recall that in the scalar case the spectrum of
each spectral band has multiplicity
2 with a possible exception in the end points of the bands. For $\cL$ this is similar.

\no iii) In the case  $V=0$ the
corresponding functions have the forms
\[
\lb{0}  \D_1^0(\l)=\cosh z, \ \ \ \D_2^0(\l)=\cos z, \ \ \ \ \ \ 
m=1,2,\ \ \ \ \ \ z=\l^{1\/4}.
\]
 Thus the function
$\D^0(\l)=\sum_0^\iy(-1)^n{(\sqrt\l)^n\/(2n)!}$ is analytic on the
two-sheeted Riemann surface $\cR^0$ of the function $\sqrt\l$, where $\sqrt 1=1$ on the first sheet $\mR_1^0$. 
We have only one resonance gap $(-\iy,0)$, since we have only  one branch point, which equals zero. We also have
\[
\lb{D0} D_\pm^0=(\cosh z\mp 1)(\cos z \mp 1),\ \ {\l_n^\pm}^0 =(\pi n)^4,\ \
{r_n^\pm }^0=-4(\pi n)^4,\ \  n\ge 1,\ \  {r_0^-}^0=0,\ {\l_0^+}^0=0.
\]
iv)  We describe the difference between the Lyapunov functions
for the operators $\cL$ and the operator $\cL_{sys}=-{d^2\/dt^2}+\cV$ on the real line where $\cV$  is a 1-periodic $2\ts 2$ matrix potential from [BBK]:

1) For the operator $\cL$ the Lyapunov function $\D_1$
is increasing  and $\D_2$ is bounded on the real line at high energy. It creates some problem to determine the asymptotics
of the spectral data for $\cL$.  Moreover, this implies that the spectrum of $\cL$ has multiplicity 2  at high energy.
For the operator $\cL_{sys}$ all Lyapunov functions are bounded on the real line.

2) The resonances for $\cL$ go to 
$-\iy$ and  the resonances for the operator $\cL_{sys}$
go to $+\iy$.

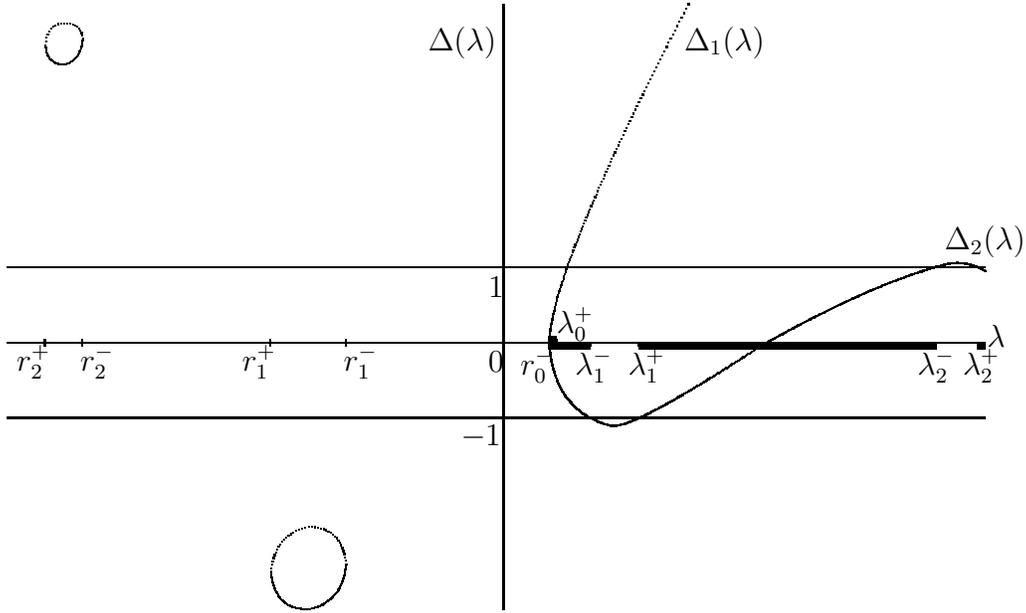
\begin{figure}
\unitlength 1.00mm \linethickness{0.4pt}
\begin{picture}(140.00,86.00)
\put(76.00,5.00){\line(0,1){80.50}}
\put(10.00,50.50){\line(1,0){130.00}}
\put(10.00,40.50){\line(1,0){130.00}}
\put(10.00,30.50){\line(1,0){130.00}}
\put(141.50,41.50){\makebox(0,0)[cc]{$\l$}}
\put(70.50,80.50){\makebox(0,0)[cc]{$\D(\l)$}}
\put(75.00,38.00){\makebox(0,0)[cc]{$0$}}
\put(73.00,28.00){\makebox(0,0)[cc]{$-1$}}
\put(75.00,48.00){\makebox(0,0)[cc]{$1$}}
\bezier{120}(82.00,40.50)(82.25,50.50)(100.50,85.50)
\put(100.00,80.50){\makebox(0,0)[lc]{$\D_1(\l)$}}
\bezier{200}(82.00,40.50)(82.50,31.50)(90.00,29.50)
\bezier{200}(90.00,29.50)(92.50,28.50)(110.00,40.00)
\bezier{200}(110.00,40.00)(125.00,48.50)(135.00,51.00)
\bezier{50}(135.00,51.00)(137.50,51.50)(140.00,50.00)
\put(140.00,54.00){\makebox(0,0)[cc]{$\D_2(\l)$}}
\bezier{40}(45.00,10.00)(45.50,5.00)(50.00,5.00)
\bezier{40}(50.00,5.00)(54.50,6.00)(55.00,11.00)
\bezier{20}(45.00,10.00)(45.50,15.00)(50.00,16.00)
\bezier{20}(50.00,16.00)(54.50,16.00)(55.00,11.00)
\put(45.00,40.00){\line(0,1){1.00}}
\put(43.50,38.00){\makebox(0,0)[cc]{$r_1^+$}}
\put(55.00,40.00){\line(0,1){1.00}}
\put(57.00,38.00){\makebox(0,0)[cc]{$r_1^-$}}
\bezier{20}(15.00,80.00)(15.25,77.50)(17.50,77.50)
\bezier{20}(17.50,77.50)(19.75,78.00)(20.00,81.00)
\bezier{10}(15.00,80.00)(15.25,83.00)(17.50,83.00)
\bezier{10}(17.50,83.00)(19.75,83.00)(20.00,81.00)
\put(15.00,40.00){\line(0,1){1.00}}
\put(13.50,38.00){\makebox(0,0)[cc]{$r_2^+$}}
\put(20.00,40.00){\line(0,1){1.00}}
\put(22.00,38.00){\makebox(0,0)[cc]{$r_2^-$}}
\put(82.00,40.00){\line(0,1){1.00}}
\put(80.50,37.50){\makebox(0,0)[cc]{$r_0^-$}}
\put(83.00,40.50){\line(0,1){0.50}}
\put(85.50,43.00){\makebox(0,0)[cc]{$\l_0^+$}}
\put(87.50,40.00){\line(0,1){0.50}}
\put(88.00,37.50){\makebox(0,0)[cc]{$\l_1^-$}}
\put(94.00,40.00){\line(0,1){0.50}}
\put(95.00,37.50){\makebox(0,0)[cc]{$\l_1^+$}}
\put(133.50,40.00){\line(0,1){0.50}}
\put(133.50,37.50){\makebox(0,0)[cc]{$\l_2^-$}}
\put(139.00,40.00){\line(0,1){0.50}}
\put(139.50,37.50){\makebox(0,0)[cc]{$\l_2^+$}}
\put(82.00,41.00){\linethickness{2.0pt}\line(1,0){1.00}}
\put(82.00,40.00){\linethickness{2.0pt}\line(1,0){5.50}}
\put(94.00,40.00){\linethickness{2.0pt}\line(1,0){39.50}}
\put(139.00,40.00){\linethickness{2.0pt}\line(1,0){1.00}}
\end{picture}
\caption{The function $\D$ for small $V$.} \lb{sp}
\end{figure}

We formulate our theorem about the asymptotics of the periodic and anti-periodic eigenvalues and resonances at high energy and
the recovering the spectrum of $\cL$.

\begin{theorem}   \lb{T2}
 Let $V\in L_0^1(\T)$. Then

\no i) There exists an integer $N\ge 0$ such that for all $n\ge N$ the
inequalities are fulfilled:
\[
\lb{T2-1}
\l_n^-\le\l_n^+<\l_{n+1}^-\le\l_{n+1}^+<\l_{n+2}^-\le\l_{n+2}^+<...
\]
where the intervals $[\l_n^+,\l_{n+1}^-]$ are spectral bands of multiplicity $2$
in  $(\l_n^+,\l_{n+1}^-)$, and the intervals $(\l_n^-,\l_n^+)$ are gaps. Moreover, the following asymptotics are fulfilled:
\[
\lb{T2-2}
\l_n^{\pm}=(\pi n)^4\pm |\hat V_n|+O(n^{-3/2}),\ \ \ \ \hat V_n=\int_0^1 V(t)e^{-i2\pi nt}dt,
\]
\[
\lb{T2-3}
r_n^{\pm}=-4(\pi n)^4\pm \sqrt{2}|\hat V_n|+O(n^{-3/2}),
\ \ \ n\to+\iy.
\]
\no ii) The periodic spectrum and the anti-periodic spectrum recover
the resonances and the spectrum of the operator $\cL$.

\no iii) The periodic (anti-periodic) spectrum is recovered by the
anti-periodic (periodic) spectrum and the resonances.
\end{theorem}

{\bf Remark.} Assume that in Theorem 1.2 the potential
$V$ has the Fourier coefficients $\hat V_n=1/n, |n|>n_*$ for some
$n_*\in \N$. Using asymptotics \er{T2-2} we deduce that
there exist infinitely many gaps in the spectrum of $\cL$ and 
infinitely many resonance gaps.
Unfortunately, we can not construct a potential
with a finite number of gaps in the spectrum of $\cL$.

Consider the operator $\cL={d^4\/dt^4}+\g V, V\in L_0^1(\T)$
and real $\g$.
We will show that for small $\g\ne 0$ the lowest spectral band of $\cL$ contains an interval $(r_0^-,\l_0^+)$ of multiplicity 4 (see Fig.\ref{sp}).

\begin{theorem} \lb{1.3}
Let $\cL={d^4\/dt^4}+\g V$, where
$V\in L_0^1(\T),V\ne 0,\g\in\R$. Then
there exist two real analytic functions $r_0^-(\g), \l_0^+(\g)$
in the disk $\{|\g|<\ve\}$ for some $\ve>0$ such that $r_0^-(\g)<\l_0^+(\g)$ for all $\g\in (-\ve,\ve)\sm\{0\}$. Here
 $r_0^{-}(\g)$ is a simple zero of the function $\r(\l,\g V), r_0^-(0)=0$ and
$\l_0^+(\g)$ is a simple zero of the function $D_+(\l,\g V), \l_0^+(0)=0$.
Moreover,  the following asymptotics are fulfilled:
\[
\lb{r0}\lb{l0}
r_0^-(\g)=2\g^2 (4v_1-v_2)+O(\g^3),\ \ \
 \l_0^+(\g)=2\g^2 (4v_1-v_2)+O(\g^3),
\]
\[
\lb{bg0} \l_0^+(\g)-r_0^-(\g)=4A^2\g^4+O(\g^5),\ \ \
A={v_2\/12}-{4\/3}v_1={5\/4}\sum_{n\ne 0}{|\hat V_n|^2\/(2\pi n)^6}>0,
\]
as $\g\to 0$, where
\[
\lb{vm} \hat V_n=\int_0^1 V(t)e^{-i2\pi n}dt,\ \ 
v_m={1\/144}\int_0^m dt\int_0^tv(s)v(t)(m-t+s)^3(t-s)^3ds,\ \ \ m=1,2.
\]
Furthermore, the spectral interval $(r_0^-(\g),\l_0^+(\g))$ has multiplicity 4
for any $\g\in (-\ve,\ve)\sm\{0\}$.
\end{theorem}

Consider the operator $\cL^\g={d^4\/dt^4} +\g \d_{per},\g\in\R,$ where $\d_{per}(t)=\sum \d(t-n)$. We prove that the function $\D(\l,\g \d_{per})$ has real and as well as non-real branch points
for some $\g>0$.

\begin{proposition} \lb{1.4}
There exists $N>0$ such that for each $n\ge N$
there exist $z_n\in(2n\pi,(2n+1)\pi),\g_n\in\R,\ve_n>0,$
and the functions $r_{n}^\pm(\g),-\ve_n<\g-\g_n<\ve_n$,
such that $r_{n}^\pm(\g)$ are zeros of the function
$\r(\l,\g \d_{per}),r_{n}^\pm(\g_n)=z_n^4$. Moreover, the following asymptotics are fulfilled:
\[
\lb{dr} r_{n}^\pm(\g)=z_n^4
\pm \a_n z_n^3\sqrt\n
+O(\n^{{3\/2}}),\ \ \ \a_n>0,\ \ \  \n=\g-\g_n\to 0.
\]
\end{proposition}

\no {\bf Remark.} i) Numerical experiments show that asymptotics \er{dr} hold for all $n\ge 1$.
The qualitative behavior of $\r(\l,\g \d_{per}),\D(\l,\g \d_{per})$ at $\g\approx \g_1$ is shown by Fig. \ref{esp}. ii) If $\n>0$, then the branch points $r_{n}^\pm(\g)$ are real. If $\n<0$, then the branch points $r_{n}^\pm(\g)$ are non-real.

\begin{figure}
\unitlength 1.00mm \linethickness{0.4pt}
\begin{picture}(155.00,140.00)
\put(7.00,103.50){\line(0,1){35.00}}
\put(5.00,120.00){\line(1,0){70.00}}
\put(68.00,117.00){\makebox(0,0)[cc]{$r_0^+$}}
\put(48.00,117.00){\makebox(0,0)[cc]{$r_1^-$}}
\put(18.00,117.00){\makebox(0,0)[cc]{$r_1^+$}}
\put(75.00,118.00){\makebox(0,0)[cc]{$\l$}}
\put(14.00,137.00){\makebox(0,0)[cc]{$\r(\l,\g \d_{per})$}}
\bezier{200}(10.00,103.50)(18.00,134.00)(26.00,133.50)
\bezier{100}(26.00,133.50)(32.00,134.00)(47.00,122.00)
\bezier{50}(47.00,122.00)(54.00,116.00)(60.00,116.00)
\bezier{200}(60.00,116.00)(66.00,116.00)(75.00,138.00)
\put(7.00,-1.00){\line(0,1){92.00}}
\put(5.00,60.00){\line(1,0){70.00}}
\put(77.00,59.00){\makebox(0,0)[cc]{$\l$}}
\put(14.00,92.00){\makebox(0,0)[cc]{$\D(\l,\g \d_{per})$}}
\put(5.00,63.00){\line(1,0){70.00}}
\put(3.00,63.00){\makebox(0,0)[cc]{$1$}}
\put(5.00,57.00){\line(1,0){70.00}}
\put(2.00,57.00){\makebox(0,0)[cc]{$-1$}}
\put(65.30,59.00){\line(0,1){2.00}}
\put(65.00,66.00){\makebox(0,0)[cc]{$r_0^+$}}
\put(49.50,59.00){\line(0,1){2.00}}
\put(49.50,66.00){\makebox(0,0)[cc]{$r_1^-$}}
\put(15.00,59.00){\line(0,1){2.00}}
\put(15.00,66.00){\makebox(0,0)[cc]{$r_1^+$}}
\bezier{200}(65.50,58.00)(66.00,65.00)(75.00,90.00)
\bezier{200}(65.50,58.00)(66.00,54.00)(75.00,56.00)
\bezier{200}(15.00,15.00)(16.00,25.00)(35.00,40.00)
\bezier{200}(35.00,40.00)(49.00,50.00)(49.50,45.00)
\bezier{200}(15.00,15.00)(16.00,7.00)(35.00,22.00)
\bezier{200}(35.00,22.00)(49.50,35.00)(49.50,45.00)
\put(50.00,0.00){\makebox(0,0)[cc]{$\g>\g_1$}}
\put(87.00,103.50){\line(0,1){35.00}}
\put(85.00,120.00){\line(1,0){70.00}}
\put(98.00,117.00){\makebox(0,0)[cc]{$r_1^+$}}
\put(155.00,118.00){\makebox(0,0)[cc]{$\l$}}
\put(94.00,137.00){\makebox(0,0)[cc]{$\r(\l,\g \d_{per})$}}
\bezier{200}(90.00,103.50)(98.00,140.00)(106.00,139.50)
\bezier{100}(106.00,139.50)(112.00,140.00)(127.00,128.00)
\bezier{50}(127.00,128.00)(134.00,121.00)(140.00,123.00)
\bezier{200}(140.00,123.00)(146.00,125.00)(155.00,140.00)
\put(87.00,-1.00){\line(0,1){92.00}}
\put(85.00,60.00){\line(1,0){70.00}}
\put(157.00,59.00){\makebox(0,0)[cc]{$\l$}}
\put(94.00,92.00){\makebox(0,0)[cc]{$\D(\l,\g \d_{per})$}}
\put(85.00,63.00){\line(1,0){70.00}}
\put(85.00,57.00){\line(1,0){70.00}}
\bezier{200}(97.00,10.00)(115.00,8.00)(123.00,22.00)
\bezier{200}(123.00,22.00)(133.00,38.00)(148.00,50.00)
\bezier{200}(148.00,50.00)(155.00,55.00)(155.00,54.50)
\bezier{15}(94.00,13.00)(94.50,10.50)(97.00,10.00)
\bezier{200}(94.00,13.00)(95.00,30.00)(110.00,37.00)
\bezier{200}(110.00,37.00)(129.00,45.00)(140.00,58.00)
\bezier{200}(140.00,58.00)(146.00,65.00)(155.00,90.00)
\put(94.00,59.00){\line(0,1){2.00}}
\put(94.00,66.00){\makebox(0,0)[cc]{$r_1^+$}}
\put(130.00,0.00){\makebox(0,0)[cc]{$\g<\g_1$}}
\end{picture}
\caption{The qualitative behavior of the functions $\r(\l,\g \d_{per})$, $\D(\l,\g \d_{per})$ at $\g\approx \g_1$.}
\lb{esp}
\end{figure}
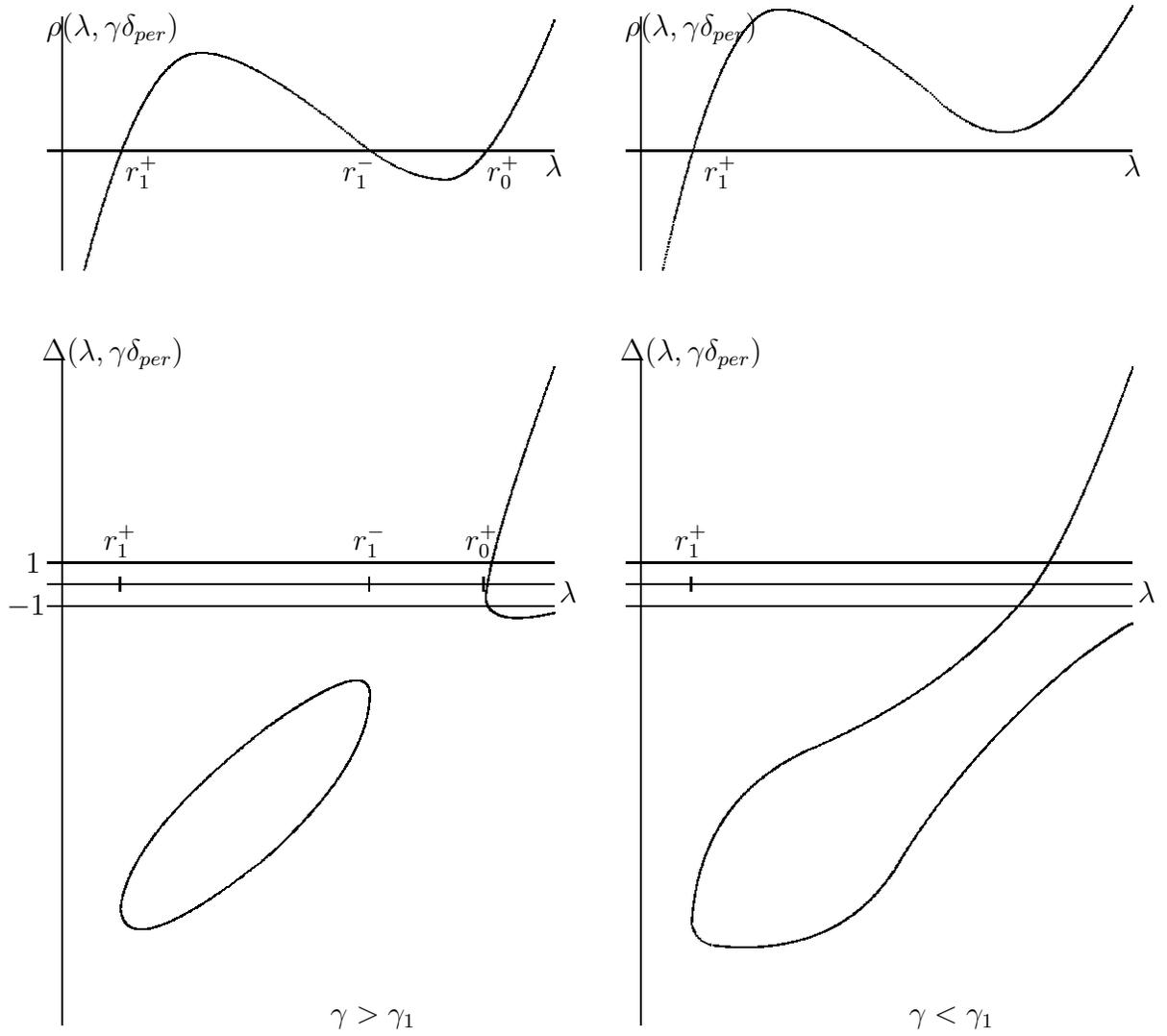

We describe the plan of our paper. In Sect.2 we obtain the basic properties of the fundamental solutions
and give a convenient representation of the functions $T_m$ in terms of the entries of
some auxiliary matrix $\F$. In Sect.3 we determine the asymptotics of the matrix $\F$. Using
these results in Sect.4 we determine the asymptotics of $T_m$. We prove Theorems \ref{T1}, \ref{T2} in Sect.5.
In Sect.6 we consider the case of small potentials and prove Theorem \ref{1.3}.
The periodic $\d$-potential is considered in  Sect.7. The existence of non-real branch points
of the Lyapunov function is shown and Proposition \ref{1.4} is proved in Section 7.

Recall that the spectrum of the operator  $-{\pa^2 \/\pa {\bf x}^2}+Q({\bf x}),\ {\bf x}\in \R^3,$ where $Q$ is a real periodic potential is absolutely continuous and consists from spectral
bands separated gaps.
We have a conjecture that  ends of these spectral bands
are periodic or anti-periodic eigenvalues or some numbers similar to resonances from our paper.

\section {Fundamental solutions}
\setcounter{equation}{0}

We begin with some notational convention. A vector
$h=\{h_n\}_1^N\in \C^N$ has the Euclidean norm
$|h|^2=\sum_1^N|h_n|^2$, while a $N\ts N$ matrix $A$ has the
operator norm given by $|A|=\sup_{|h|=1} |Ah|$. In this section
we study the fundamental solutions $\vp_j,j=0,1,2,3,$. We introduce the
fundamental solutions $\vp_j^0$ of the unperturbed equation
$y''''=\l y$ given by
\[
\label{2.1a} \vp_0^0(t,\l)={\cosh zt+\cos zt\/2},\ \ \
\vp_1^0(t,\l)={\sinh zt+\sin zt\/2z},
\]
\[
\label{6} \vp_2^0(t,\l)={\cosh zt-\cos zt\/ 2z^2},\ \ \
\vp_{3}^0(t,\l)={\sinh zt-\sin zt\/2z^3},\ \ \ z=\l^{1/4}, \
\arg z\in (-{\pi\/4},{\pi\/4}],
\]
which are entire in $\l\in \C$. Here below we have $z=x+iy, \ x\ge|y|$. They satisfy
\[
\label{(2.4)}
\pa_t^k\vp_{j}^0(t,\l)=\vp_{j-k}^0(t,\l),\ \ \
\sum_{m=0}^3\vp_{j-m}^0(t,\l)\vp_{m-k}^0(s,\l)=\vp_{j-k}^0(t+s,\l),
\ \ \ 0\le k,j\le 3,
\]
\[
\lb{itn} \sum_{j=0}^3\vp_{0,j}^{(j)}(m,\l)=2(\cosh mz+\cos mz),\
\ \ \ \  m\ge 1,
\]
\[
\lb{ep} |\vp_{j}^0(t,\l)|\le {e^{tx}+e^{t|y|}\/2}\le e^{xt}
,\ \ \ j=0,1,2,3.
\]
The fundamental solutions $\vp_j,j=0,1,2,3,$ satisfy the following integral equations
\[
\lb{2ic} \vp_j(t,\l)=
\vp_{j}^0(t,\l)-\int_0^t\vp_{3}^0(t-s,\l)V(s)\vp_j(s,\l)ds,\ \ \
\ \ (t,\l)\in\R\ts\C.
\]
The standard iterations in \er{2ic} yield
\[
\lb{2id} \vp_j(t,\l)=\sum_{n\ge 0}
\vp_{n,j}(t,\l),\ \ \ \vp_{n+1,j}(t,\l)=
-\int_0^t\vp_{3}^0(t-s,\l)V(s)\vp_{n,j}(s,\l)ds,
\]
where $\vp_{0,j}=\vp_{j}^0$.  Define the functions
\[
\label{(7.1)}
T_{m,2}(\l)={1\/4}\int_0^m
dt\int_0^tV(s)V(t)\vp_{3}^0(m-t+s,\l)\vp_{3}^0(t-s,\l)ds,\ \ \
m=1,2.
\]
We prove

\begin{lemma} \lb{T21}
 For each $(t,V)\in\R_+\ts L_0^1(\T)$ and $j=0,1,2,3$ the functions $\vp_j(t,\cdot)$ are real on $\R$ and entire
and for each $N\ge -1$ the following estimates are fulfilled:
\[
\lb{2if}
\max_{0\le j,k\le 3}\lt\{\lt|\l^{j-k\/4}\lt(\vp_j^{(k)}(t,\l)-
\sum_0^N\vp_{n,j}^{(k)}(t,\l)\rt)\rt|\rt\} \le
{(\vk t)^{N+1}\/(N+1)!}e^{xt+\vk},\ \ \ \ \ \ \vk={\|V\|\/{|\l|_1^{3/4}}},
\]
where $|\l|_1\ev\max\{1,|\l|\}$.
Moreover, $T_m, m=1,2$ is real for $\l\in\R$, entire and
satisfies
\[
\label{2.2a}
|T_m(\l)|\le e^{xm+\vk},\ \
 |T_m(\l)-T_{m}^0(\l)|\le {(m\vk)^2\/2}e^{xm+\vk},
\]
\[
|T_m(\l)-T_{m}^0(\l)-T_{m,2}(\l)|\le {(m\vk)^3\/3!}e^{xm+\vk}.
\label{(4.61)}
\]
\end{lemma}
\no {\it Proof.}  We estimate $\vp_0$, the proof of
other estimates is similar. \er{2id} gives
\[
\vp_{n,0}(t,\l)= \int\limits_{0< t_n<...< t_2< t_1\le t_0=t}
\lt(\prod\limits_{1\le k\le n}
\vp_{3}^0(t_{k-1}-t_k,\l)V(t_k)\rt)\vp_{0,0}(t_n,\l)dt_1dt_2...dt_n.
\lb{2ig}
\]
Substituting estimates \er{ep} into \er{2ig}
we obtain $| \vp_{n,0}(t,\l)|\le
{(\vk t)^n\/n!}e^{xt},$ which shows that for any fixed
$t\in[0,1]$ the formal series \er{2id} converges uniformly on
bounded subset of $\C$. Each term of this series is an entire
function. Hence the sum is an entire function. Summing the
majorants we obtain estimates  \er{2if}.

The monodromy matrix is real on the real line. Then
$T_1,T_2$ are real on $\R$. We will prove \er{2.2a}, \er{(4.61)}.
We have
\[
\lb{mum}
4T_m=\Tr M^m(\l)=\Tr M(m,\l)=\sum_{j=0}^3\vp_j^{(j)}(m,\l)
=\sum_{n\ge 0}\sum_{j=0}^3\vp_{n,j}^{(j)}(m,\l),\ \ \ m=1,2.
\]
The estimates $|\vp_{n,j}^{(j)}(m,\l)|\le {(m\vk)^n\/n!}e^{xm+\vk}$ yield
\[
\lb{svp}
\lt|\sum_{j=0}^3\vp_{n,j}^{(j)}(m,\l)\rt|\le 4{(m\vk)^n\/n!}e^{xm+\vk},\ \ n\ge 0.
\]
The last estimate shows that the series  \er{mum} converges
uniformly on bounded subset of $\C$. Each term of this series is
an entire function. Hence the sum is an entire function and
$T_1,T_2$ are entire. Summing the majorants we obtain the first
estimates in \er{2.2a}.
Using \er{(2.4)}, \er{2id} we obtain
$$
\sum_{j=0}^3\vp_{1,j}^{(j)}(m,\l)
=-\sum_{j=0}^3\int_0^m\vp_{3-j}^0(m-t,\l)\vp_{j}^0(t,\l)V(t)dt
=-\vp_{3}^0(m,\l)\int_0^mV(t)dt=0,
$$
and
$$
\sum_{j=0}^3\vp_{2,j}^{(j)}(m,\l)=\sum_{j=0}^3\int_0^m\int_0^t
\vp_{3-j}^0(m-t,\l)\vp_{3}^0(t-s,\l)\vp_{j}^0(s,\l)V(s)V(t)dsdt
=T_{m,2}(\l),
$$
where we have used \er{(7.1)}. Then \er{mum}, \er{svp} give the
second estimate in \er{2.2a} and \er{(4.61)}.$\BBox$

Note that $\mM(t,\l)$ is a solution of the equation
\[
\lb{Me} Y'=\left(\begin{array}{cccc}
0&1&0&0\\
0&0&1&0\\
0&0&0&1\\
\l-V&0&0&0\\
\end{array}\right)Y,\ \ \ \ \  \ (t,\l)\in[0,1]\ts\C,
\]
such that $\mM(0,\l)=I_4$. Lemma 2.1 does not give
asymptotics of $\mM(1,\l)$. In order to determine the asymptotics 
of $\mM(1,\l)$ we need another solution $Y(t,\l)$ of Eq.\er{Me}
with the good asymptotics at high energy, see Lemma 2.2. Note that
$\mM(t,\l)=Y(t,\l)Y^{-1}(0,\l)$.

We will construct the matrix $Y$ using some special solutions
of Eq.\er{1b}, which have "good" asymptotics at $|\l|\to\iy$ (see Lemma \ref{T31}).

In order to determine the asymptotics 
of $M(\l)$ we need another solution $\vt_j,j=0,1,2,3$.
We introduce a matrix $\O(\l)$ by 
$$
\O=\O(\l)=\diag(\o_0,\o_1,\o_2,\o_3)=(1,-i,i,-1),\ \ \ \l\in\ol{\C_+},\ \ \
\O(\ol\l)=\ol{\O(\l)}.
$$
Here $\o_0$ and $\o_3$ are constants in $\C$, but $\o_j=\o_j(\l),j=1,2,$ are constants only in $\C_{\pm}$.
It will imply that some of functions, which will be introduced below, they will not be analytic in whole complex plane, but will be analytic only in $\C_\pm$ or theirs subsets.

We define the functions $a_k(t,\l)$, $k=0,1,2,3, (t,\l)\in\R\ts\ol{\C_\pm}$ by
\[
\lb{dkm} a_0(t,\l)=0,\ \ a_k(t,\l)=\sum_{j=0}^{k-1}\o_je^{zt\o_j},\ \ k=1,2,3, \ \ \ t<0,
\]
\[
\lb{dk} a_k(t,\l)=-\sum_{j=k}^{3}\o_je^{z t\o_j}, \ \ \
k=0,1,2,3,\ \ \ t>0,
\]
$z=\l^{1/4},-{\pi\/4}<\arg z<{\pi\/4}$.
Note that $\Re(z\o_0)\ge\Re(z\o_1)\ge\Re(z\o_2)\ge\Re(z\o_3)$ (see Fig.3).
\begin{figure}
\unitlength 1.00mm \linethickness{0.4pt}
\begin{picture}(81.00,77.00)
\put(77.00,8.00){\line(0,1){69.00}}
\put(43.00,41.00){\line(1,0){68.00}}
\put(77.00,41.00){\line(-1,-1){30.00}}
\put(77.00,41.00){\line(1,1){33.00}}
\put(77.00,41.00){\line(1,-1){30.00}}
\put(77.00,41.00){\line(-1,1){32.00}}
\put(110.00,39.00){\makebox(0,0)[cc]{$x$}}
\put(75.00,75.00){\makebox(0,0)[cc]{$y$}}
\put(107.00,44.00){\makebox(0,0)[cc]{$z\o_0$}}
\put(107.00,47.00){\rule{1.00\unitlength}{1.00\unitlength}}
\put(67.00,71.00){\makebox(0,0)[cc]{$z\o_2$}}
\put(71.00,71.00){\rule{1.00\unitlength}{1.00\unitlength}}
\put(47.00,32.00){\makebox(0,0)[cc]{$z\o_3$}}
\put(47.00,35.00){\rule{1.00\unitlength}{1.00\unitlength}}
\put(89.00,11.00){\makebox(0,0)[cc]{$z\o_1$}}
\put(83.00,11.00){\rule{1.00\unitlength}{1.00\unitlength}}
\put(77.00,2.00){\makebox(0,0)[cc]{$z=x+iy$-plane}}
\end{picture}
\caption{}
\end{figure}
Then we have
$$
|e^{z t\o_j}|=e^{t\Re(z\o_j)}\le e^{t\Re(z\o_k)},\ \ \ 0\le j\le k, \ \ \ t<0,\ \ \
|e^{z t\o_j}|\le e^{t\Re(z\o_k)},\ \ \ k\le j\le 3,\ \ \ t\ge 0.
$$
Substituting these estimates into identities \er{dkm},\er{dk} we obtain
\[
\lb{ek} |a_k(t,\l)|\le 4{e^{t\Re(z\o_k)}},\ \ \ k=0,1,2,3, \ \ \ \
(t,\l)\in \R\ts\ol{\C_\pm}.
\]
Identities \er{dkm}, \er{dk} show that $a_k,a_k',a_k''$ are
continuous functions of $t\in \R$ and
\[
\lb{pak} a_k'''(+0,\l)-a_k'''(-0,\l)=-4z^3,\ \ \ a_k''''(t,\l)-\l
a_k(t,\l)=-4z^3\d(t).
\]
Let $\L_r=\{\l\in\C:r\|V\|<|\l|^{3/4}\}$ and $\L_r^{\pm}=\L_r\cap\C^{\pm}, r>0$. Below we need

\begin{lemma} \lb{T31}
Let $V\in L_0^1(\T)$. Then for each $j=0,1,2,3$ and $\l\in\ol{\L_1^\pm}$
the integral equation
\[
\lb{iey}
\vt_j(t,\l)=
e^{zt\o_j}+{1\/4z^3}\int_0^1a_j(t-s,\l)V(s)\vt_j(s,\l)ds,\ \ \
t\in [0,1],
\]
has the unique solution $\vt_j(\cdot,\l)$ and each $\vt_j(t,\cdot),
t\in [0,1]$ is analytic in $\L_1^\pm$ and continuous in $\ol{\L_1^\pm}$. Moreover, each
$\vt_j(\cdot,\l),\l\in\ol{\L_1^\pm}$ is a solution of equation
$\vt_j''''+V\vt_j=\l \vt_j$ for $t\in [0,1]$  and satisfies
\[
\lb{ee1}
 |\vt_j(t,\l)|\le {e^{t\Re(z\o_j)}\/1-\vk},\ \ \
 (t,\l)\in [0,1]\ts\ol{\L_1^\pm},\ \ \ \vk={\|V\|\/|z|^3}.
\]
\end{lemma}

\no {\it Proof.}
Let $\vt_{j,0}(t,\l)=e^{zt\o_j}$.
The iterations in Eq.\er{iey} provide the identities
 \[
\lb{ef}\lb{ynj} \vt_j(t,\l)=\sum_0^\iy \vt_{j,n}(t,\l),\ \ \
\vt_{j,n}(t,\l)={1\/4z^3}\int_0^1a_j(t-s,\l)V(s)\vt_{j,n-1}(s,\l)ds,\
\ n\ge 1,
\]
$$
\vt_{j,n}(t,\l)={1\/(4z^3)^n}\int_{[0,1]^n}a_j(t-t_n)a_j(t_n-t_{n-1})...
a_j(t_2-t_1)e^{zt\o_j}V(t_1)...V(t_n)dt_1...dt_n.
$$
For $-{\pi\/4}<\arg z<{\pi\/4}, t\in [0,1]$ using estimates \er{ek} we obtain
$$
|\vt_{j,n}(t,\l)|\le {1\/|z|^{3n}}\int\limits_{[0,1]^n}
e^{(t-t_n)\Re(z\o_j)}e^{(t_n-t_{n-1})\Re(z\o_j)}...
e^{(t_2-t_1)\Re(z\o_j)}e^{t_1\Re(z\o_j)}
$$
$$
\ts|V(t_1)|...|V(t_n)|dt_1...dt_n
={e^{t\Re(z\o_j)}\/|z|^{3n}}\int_{[0,1]^n}
|V(t_1)|...|V(t_n)|dt_1...dt_n={\|V\|^n\/|z|^{3n}}e^{t\Re(z\o_j)}.
$$
Substituting the last estimate
into \er{ef} we get \er{ee1}. This estimate shows that for each
fixed $|z|>\|V\|_1^{1/3},-{\pi\/4}<\arg z<{\pi\/4}$ series
\er{ef} converges uniformly on the interval $[0,1]$. Thus it
gives the solution of Eq.\er{iey}. Suppose that there exists
another solution $\wt y_j$ of this equation. Then
$y=y_j-\wt y_j$ satisfies the equation
$y(t,\l)={1\/4z^3}\int_0^1V(s)a_j(t-s,\l)y(s,\l)ds$. The
iterations of this equation give $y(t,\l)=0$. For each $t\in\R$
the series \er{ef} converges uniformly on any bounded subset of
the domain $|z|>\|V\|^{1/3},-{\pi\/4}<\arg z<0$, or $0<\arg z<{\pi\/4}$. Each term of this
series is an analytic function of $z$ in this domain. Hence the
sum is an analytic function and $\vt_n(t,\l)$ is analytic  of
$\l \in\L_1^\pm$ and is continuous in $\ol{\L_1^\pm}$.

We will show that $\vt_j(\cdot,\l)$ is a solution of $y''''+Vy=\l y$.
Using \er{pak} and \er{ynj} we obtain
$$
\vt_j''''(t,\l)=\l\vt_j^0(t,\l)+\!\!\int_0^1\lt({\l
a_j(t-s,\l)\/4z^3}-\d(t-s)\rt)V(s)\vt_{j}(s,\l)ds=
\l\vt_j(t,\l)-V(t)\vt_{j}(t,\l)
$$
Thus $\vt_j(\cdot,\l),\l\in\ol{\L_1^{\pm}}$ is a solution of $y''''+Vy=\l y$, $t\in [0,1]$.
 $\BBox$

Let $Y(t,\l)=\{\vt_j^{(k)}(t,\l)\}_{k,j=0}^3$.
Then $Y(t,\l)$ is a solution of Eq.\er{Me} and hence
\[
\lb{YM} \mM(t,\l)=Y(t,\l)Y^{-1}(0,\l),\ \ \ \l\in\L_1^\pm.
\]
Differentiating Eq.\er{iey} we obtain
\[
\lb{ypr} \vt_j^{(k)}(t,\l)=(z\o_j)^ke^{zt\o_j}+
{1\/4z^3}\int_0^1V(s)a_j^{(k)}(t-s,\l)\vt_j(s,\l)ds, \ \ \
j,k=0,1,2,3.
\]
We introduce the matrix $\P=\{\p_{kj}\}_{k,j=0}^3$, where
$\p_{kj}(t,\l)=z^{-k}e^{-zt\o_j}\vt_j^{(k)}(t,\l).$
Then
\[
\lb{YP} \P(t,\l)=Z^{-1}(\l)Y(t,\l)e^{-zt\O},\ \ \ \ Z(\l)=\diag(1,z,z^2,z^3),
\]
hence $Y(t,\l)=Z(\l)\P(t,\l)e^{zt\O}$. Substituting this identity into \er{YM}, we obtain
\[
\lb{MY} \mM(t,\l)=Z(\l)\P(t,\l)e^{zt\O}\P^{-1}(0,\l)Z^{-1}(\l),\ \ \ \ \ \ t\in [0,1],
\]
and then
\[
\lb{TrM} \Tr M(\l)=\Tr(\F(\l)e^{z\O}),
\ \ \ \Tr M^2(\l)=\Tr(\F(\l) e^{z\O})^2,\ \ \ \F(\l)=\P^{-1}(0,\l)\P(1,\l),
\]
for $\l\in \ol{\L_1^\pm}$.
Note that $\Tr M$, $\Tr M^2$ are entire functions, but the functions $\P(t,\cdot),t\in\R$, and $\F$, are not entire.

\section {Asymptotics of the matrices $\P$ and $\F$ as $|\l|\to\iy$}
\setcounter{equation}{0}

We rewrite identities \er{ypr} in the form
\[
\lb{ynm}
\p_{kj}(t,\l)=\o_j^{k}+{1\/4z^3}\int_0^1a_{kj}(t-s,\l)V(s)\p_{0j}(s,\l)ds,\ \ \
a_{kj}(t,\l)=z^{-k}e^{-z\o_jt}a_j^{(k)}(t,\l),
\]
$j,k=0,1,2,3.$ The last identities contain the equations for the functions $\p_{0j}(t,\l)$:
\[
\lb{yn1}
\p_{0j}(t,\l)=1+{1\/4z^3}\int_0^1a_{0j}(t-s,\l)V(s)\p_{0j}(s,\l)ds,
\ \ \ j=0,1,2,3.
\]
We rewrite  \er{yn1} in terms of the matrix
$\p=\diag(\p_{00},\p_{01},\p_{02},\p_{03})$ by
\[
\lb{ieU}\lb{R}
\p(t,\l)=I_4+{1\/4z^3}\int_0^1a(t-s,\l)V(s)\p(s,\l)ds,\ \ \
a=\diag(a_{00},a_{01},a_{02},a_{03}).
\]
Define the matrix
$A=\{a_{kj}\}_{k,j=0}^3$ and rewrite \er{ynm} in the form
\[
\lb{ieY} \P(t,\l)=\P_0+{1\/4z^3}\int_0^1A(t-s,\l)V(s)\p(s,\l)ds,\
\ \ \ \P_0=\{{\P_0}_{kj}\}_{k,j=0}^3=\{\o_j^{k}\}_{k,j=0}^3.
\]
Note that $\P_0=\P_0(\l)$ is constant in $\ol{\C_\pm}$ and satisfies the identity
 $\P_0\P_0^*=4I_4$, hence ${1\/2}\P_0$ is an unitary matrix.
Identities \er{dkm} yield
$$
a_j^{(k)}(t,\l)=\sum_{p=0}^{j-1}\o_p(z\o_p)^ke^{zt\o_p},\ \ t<0,
\ \ \ {\rm and }\ \ \ 
a_j^{(k)}(t,\l)=-\sum_{p=j}^{3}\o_p(z\o_p)^ke^{zt\o_p},\ \ t>0.
$$
Then \er{ynm} gives
$$
a_{kj}(t,\l)=\sum_{p=0}^{j-1}\o_p\o_p^ke^{zt(\o_p-\o_j)},\ \ t<0,\ \ \ \ \ \ {\rm and }\ \ \ 
a_{kj}(t,\l)=-\sum_{p=j}^{3}\o_p\o_p^ke^{zt(\o_p-\o_j)},\ \ t>0.
$$
Thus the following identities are fulfilled:
\[
\lb{H} A(t,\l)=\P_0 \O H(t,\l),\ \ \O=\diag(\o_0,\o_1,\o_2,\o_3),
\ \ \  \ \ \ H=\{h_{jk}\}_{j,k=0}^3,
\]
where
\[
\lb{He}\lb{He1} h_{jk}(t,\l)=
\left\{\begin{array}{ccc}
0&,&j\ge k\\
e^{zt(\o_j-\o_k)}&,&j<k\\
\end{array}\right.,\ \ \ t<0,\ \ \
h_{jk}(t,\l)=
\left\{\begin{array}{ccc}
-e^{zt(\o_j-\o_k)}&,&j\ge k\\
0&,&j<k\\
\end{array}\right.,\ \ \ t\ge0.
\]
Iterations in Eq. \er{ieU} yield
\[
\lb{sU} \p(t,\l)=\sum_0^\iy \p_n(t,\l),\ \ \ \ \ \
\p_n(t,\l)={1\/4z^3}\int_0^1a(t-s,\l)V(s)\p_{n-1}(s,\l)ds, \ \ \ \p_0=I_4.
\]
Substituting the series \er{sU} into \er{ieY} we obtain
\[
\lb{sY} \P(t,\l)=\sum_0^\iy \P_n(t,\l),\ \ \ \ \ \
\P_n(t,\l)={1\/4z^3}\int_0^1A(t-s,\l)V(s)\p_{n-1}(s,\l)ds.
\]
We need the result about the matrix functions $\P, \ 
\F(\l)=\P^{-1}(0,\l)\P(1,\l),\ \l\in\L_1^\pm$.

\begin{lemma}\lb{T32} Let $V\in L_0^1(\T)$ and let $\vk={\|V\|\/|\l|^{3/4}}$. Then

\no i) Each matrix function $\P(t,\cdot),t\in [0,1]$ is analytic in
the domain $\L_2^\pm$ and continuous in $\ol{\L_2^\pm}$ and satisfies
\[
\lb{eYn} |\P_0|=2,\ \ \ |\P_n(t,\l)|\le 2\vk^n,\ \ \ n\ge 1, \ \ \
\]
\[
\lb{eY} |\P(t,\l)|\le 4,\ \ \ |\P(t,\l)-\sum_0^N
\P_n(t,\l)|\le 4\vk^{N+1},\ \ \ N\ge 0,\ \ \ \l\in\ol{\L_2^\pm}.
 \]
\no ii)  The matrix function $\F(\l)=\P^{-1}(0,\l)\P(1,\l)$ is analytic in $\l\in\L_{4}^\pm$ and continuous in $\ol{\L_4^\pm}$ and satisfies
\[
\lb{eW1}\lb{W1}
\F=I_4+\F_1+\F_2+\wt\F ,\ \
\F_1(\l)={\O\/4z^3}\int_0^1V(s)\lt(H(1-s,\l)-H(-s,\l)\rt)ds,
\]
\[
\lb{iW2}
\F_2(\l)=
{\O\/16z^6}\int_0^1du\int_0^1V(u)V(s)\lt(F(0,u,s,\l)-F(1,u,s,\l)\rt)ds,
\]
$$
F(t,u,s,\l)=H(-u,\l)\O H(t-s,\l)-H(t-u,\l)a(u-s,\l),
$$
\[
\lb{eW2}
|\F(\l)|\le 4,\ \ \
|\F_1(\l)|\le 2\vk,\ \ \ \ |\F_2(\l)|\le 4\vk^2,\ \ \ \ \
|\wt\F(\l)|\le 60\vk^3,\ \ \ \l\in\ol{\L_{4}^\pm}.
\]
\end{lemma}

\no {\it Proof.} i) Recall that ${1\/2}\P_0$ is a unitary matrix. Then $|\P_0|=2$.
Let $a(t)=a(t,\l)$. Identity \er{sU} for $\p_n$ gives
$$
\p_n(t,\l)={1\/(4z^3)^n}\int_{[0,1]^n}a(t-t_1)a(t_1-t_2)...
a(t_{n-1}-t_{n})V(t_1)...V(t_n)dt_1...dt_n, \ n\ge 1.
$$
Identity \er{ynm} for $a_{kj}$ together with estimate \er{ek} yields
$|a_{0j}(t)|\le 4$. 
Then $|a(t,\l)|\le 4$ for $(t,\l)\in \R\ts\ol{\C_\pm}$.
Using the last estimate we obtain
\[
\lb{eUn} |\p_n(t,\l)|\le \vk^n,\ \ \ n\ge 0, \ \ \
(t,\l)\in\ [0,1]\ts\ol{\C_\pm}.
\]
Identities \er{H}-\er{He1} imply
$|H(t,\l)|\le 4,|\O|=1$.
Identity \er{H} for $A$ gives $|A(t,\l)|\le 8$. Substituting
this estimate and estimates \er{eUn} into identity \er{sY} for
$\P_n$ we get \er{eYn}.
Estimate \er{eYn} shows that for each $t\in [0,1]$ the series
\er{sY} converges uniformly on any bounded domain in $\L_{1}^\pm$. Each
term in this series is an analytic function in $\L_{1}^\pm$. Hence
$\P(t,\l)$ is also analytic function in $\L_{1}^\pm$. Summing the
majorants we obtain
$$
|\sum_{N+1}^\iy\P_n(t,\l)|\le 2\sum_{N+1}^\iy\vk^n={2\vk^{N+1}\/1-\vk}\le 4\vk^{N+1},
\ \ \ \l\in\ol{\L_2^\pm},\ \ \ N\ge -1,
$$
since $\vk\le{1\/2},\l\in\ol{\L_2^\pm}$. Then we obtain \er{eY}.

\no ii) For the case $N=2$ \er{eY} yields
\[
\lb{Psi} \P(t,\l)=\P_0+\P_1(t,\l)+\P_2(t,\l)+\wt\P(t,\l),\ \ \ |\wt\P|\le 4\vk^3,
\ \ \ (t,\l)\in[0,1]\ts\ol{\L_2^\pm}.
\]
We introduce the matrices
$\P^0=\P(0,\cdot),R=(\P^0)^{-1},\P_n^0=\P_n(0,\cdot),\wt\P^0=\wt\P(0,\cdot)$.
We will prove that the matrix $R(\l)$ is analytic in the domains
$\L_{4}^\pm$ and satisfies
\[
\lb{sYi} R=R_0+R_1+R_2+\wt R,\ \ \ R_0=\P_0^{-1},\ \ \ R_1=-R_0\P_1^0R_0,\ \ \
R_2=-R_0\P_2^0R_0+R_0(\P_1^0R_0)^2,
\]
\[
\lb{wR} \wt R=-R_0\wt\P^0 R_0+R_0(\P_2^0+\wt\P^0)R_0(\P^0-\P_0)R_0
+R_0\P_1^0R_0(\P_2^0+\wt\P^0)R_0
-R((\P^0-\P_0)R_0)^3,
\]
\[
\lb{eRn}\lb{eYi}
|R|\le 1,\ \ \ |R_0|={1\/2},\ \ \ |R_1|\le{\vk\/2},\ \ \ |R_2|\le \vk^2,
\ \ \ |\wt R|\le 13\vk^3,\ \ \ \l\in\ol{\L_4^\pm}.
\]
The matrix ${1\/2}\P_0$ is unitary, then $|R_0|=|\P_0^{-1}|=|{1\/4}\P_0^*|={1\/2}$.
We have $\P^0=\P_0+(\P^0-\P_0)$. Recall $|\P^0|=2$.
If $\l\in\L_4$, then $\vk<{1\/4}$ and \er{eY} yields $|\P^0-\P_0|\le 4\vk<1$.
Hence $0$ is not an eigenvalue of $\P^0$ (see [Ka]),
$\P^0$ is invertible and $R$ is analytic in $\L_4^\pm$.
Substituting estimates \er{eYn} into \er{sYi} we have estimates of $|R_1|,|R_2|$ in \er{eRn}.
Using the standard identity $A^{-1}-B^{-1}=-A^{-1}(A-B)B^{-1}$ for the  matrices $A=\P^0,B=\P_0$ we obtain
\[
\lb{HI} R=R_0-R(\P^0-\P_0)R_0
=R_0-R_0(\P^0-\P_0)R_0+R_0\lt((\P^0-\P_0)R_0\rt)^2-R\lt((\P^0-\P_0)R_0\rt)^3
\]
which yields \er{sYi}, \er{wR}.
 Using the first identity in \er{HI} we obtain
$$
|R|\le{|R_0|\/1-|\P^0-\P_0||R_0|}\le 1,\ \ \ \l\in\ol{\L_4^\pm},
$$
which yields estimate of $|R|$ in \er{eRn}.
Substituting this estimate and  \er{eYn}, \er{eY}, \er{Psi} into \er{wR}
we obtain estimate for $|\wt R|$ in \er{eRn}. Thus relations \er{sYi}-\er{eYi} have been proved.

Now we prove \er{W1}-\er{eW2}.
Estimates \er{eY} and \er{eYi} give \er{eW2} for $\F$.
 Let $\P=\P(1,\l),\P_n=\P_n(1,\l)$ and $R=R(\l),R_n=R_n(\l)$.
Identities \er{Psi} and \er{sYi} give \er{W1}, where
$$
\F_1=R_0\P_1+R_1\P_0,\ \ \
\F_2=R_2\P_0+R_1\P_1+R_0\P_2,\ \ \
\wt\F=R_0\wt\P+R_1(\P_2+\wt\P)+R_2(\P-\P_0)+\wt R\P.
$$
Using estimates \er{eYn}, \er{eY}, \er{eRn} we get \er{eW2}.

Let $A(t)=A(t,\l),H(t)=H(t,\l), a(t)=a(t,\l)$.
Using $R_1=-R_0\P_1^0\P_0^{-1}$ and \er{sY}  we get
$$
\F_1=R_0\P_1+R_1\P_0=R_0(\P_1(1,\l)-\P_1(0,\l))=
{R_0\/4z^3}\int_0^1V(s)(A(1-s)-A(-s))ds.
$$
Substituting \er{H} for $A$ into the last identity
we obtain \er{W1}.

We will prove \er{iW2}. Recall $\F_2=R_2\P_0+R_1\P_1+R_0\P_2$.
Identities \er{sYi} give
$R_2=((R_0\P_1^0)^2-R_0\P_2^0)R_0.$ Then
$$
R_2\P_0=(R_0\P_1^0)^2-R_0\P_2^0=R_0(\P_1^0R_0\P_1^0
-\P_2^0).
$$
Substituting $\P_1,\P_2$ from \er{sY} and  $\p_0,\p_1$ from
\er{sU} into the last identity we obtain
$$
R_2\P_0={R_0\/16z^6}\!\!\iint_{[0,1]^2}\!\!\!
V(u)V(s)A(-u)\lt(R_0A(-s)-
a(u-s)\rt)duds.
$$
Using identity \er{H} for $A$ we have
$$
 R_2\P_0={\O\/16z^6}\!\!\iint_{[0,1]^2}\!\!\!
V(u)V(s)\lt(H(-u)\O H(-s) -H(-u)a(u-s)\rt)duds.
$$
Recall that $R_1=-R_0\P_1^0R_0$. Then \er{sY} implies
$$
R_1\P_1\!=\!{-R_0\/16z^6}\!\!\iint\limits_{[0,1]^2}\!\!\!
V(u)V(s)A(-u)R_0A(1-s)duds
\!=\!{-\O\/16z^6}\!\!\iint\limits_{[0,1]^2}\!\!\!V(u)V(s)H(-u)\O H(1-s)duds.
$$
\er{sY} gives
$$
R_0\P_2={R_0\/16z^6}\!\!\iint_{[0,1]^2}\!\!\!
V(u)V(s)A(1-u)a(u-s)duds
={\O\/16z^6}\!\!\iint_{[0,1]^2}\!\!\!\!\!V(u)V(s)H(1-u)a(u-s)duds.
$$
The last three identities yield \er{iW2}.
$\BBox$

\section{Asymptotics of the trace of the monodromy matrix}
\setcounter{equation}{0}

We introduce the functions $b_{jk},c_{jk}, \a_{j}, \b_{jk}$ in  $\ol{\C_\pm}$ by
\[
\lb{bnk}
b_{jk}={\o_j\o_k\/16z^6}
\int_0^1\!\!du\int_0^u\!\!V(u)V(s)e^{z(u-s)(\o_j-\o_k)}ds,
\ \ \ c_{jk}=e^{z(\o_j-\o_k)}(b_{jk}+b_{kj}),
\]
\[
\lb{al0}\lb{al3}\lb{aln} \a_{0}=
\sum_{k=1}^3(b_{k0}+c_{k0}),\ \ \a_j=\sum_{k=j+1}^3
(b_{kj}+c_{kj})
-\sum_{k=0}^{j-1}b_{jk},\ \ j=1,2,
\ \ \ \a_{3}=- \sum_{k=0}^{2}b_{3k},
\]
\[
\lb{bt}
\b_{jk}=\a_{j}+\a_{k}-c_{kj}.
\]
which are analytic in $\C_\pm$  and are continuous in $\ol\C_\pm$. We define the function 
$$
T=4T_1^2-T_2, \ \ \ \ \ \ \ \ 
T^0=4(T_1^0)^2-T_2^0=1+2\cosh z\cos z.
$$
In order to determine asymptotics of $T_1,T$,
we will show the following identities
\[
\lb{T1T} T_1={1\/4}\sum_0^3\f_{kk}e^{z\o_k},\ \ \
T={1\/2}\sum_{0\le j<k\le
3}v_{jk}e^{z(\o_j+\o_k)},\ \ \
v_{jk}=\f_{jj}\f_{kk}-\f_{jk}\f_{kj},
\]
for $\l\in\ol{\L_4^\pm}$ where $\F=\{\f_{jk}\}_{j,k=0}^3$.  Identity  \er{TrM} yields the first identity in \er{T1T}.
Due to \er{TrM} we have
$$
T_2={1\/4}\Tr M^2={1\/4}\Tr(\F e^{z\O})^2={1\/4}\lt(\sum_0^3\f_{jj}^2e^{2z\o_j}
+2\sum_{0\le j<k\le 3}\f_{jk}\f_{kj}e^{z(\o_j+\o_k)}\rt).
$$
Substituting the last identity and the first identity in \er{T1T} into $T=4T_1^2-T_2$, we obtain the second identity in \er{T1T}. Recall $z=\l^{1/4}=x+iy,|y|\le x,\l\in\C.$

\begin{lemma}\lb{T33}
Let $V\in L_0^1(\T)$ and let $\F=\P^{-1}(0,\cdot)\P(1,\cdot)=\{\f_{jk}\}_{j,k=0}^3$. Then

\no i) For each $\l\in\ol{\L_4^\pm}$ the following identities and estimates are fulfilled
\[
\lb{efj} \f_{kk}(\l)=1+\a_k(\l)+\wt \f_{kk}(\l),\ \ \ |\a_k(\l)|\le{3\/8}\vk^2,
\ \ \ |\wt \f_{kk}(\l)|\le 120\vk^3,\ \ \ k=0,1,2,3,
\]
\[
\lb{vjk} v_{jk}(\l)=1+\b_{jk}(\l)+\wt v_{jk}(\l),\ \ \ |\b_{jk}(\l)|\le{7\/8}\vk^2,
\ \ \ |\wt v_{jk}(\l)|\le 1701\vk^3,\ \ \ 0\le j<k\le 3.
\]

\no ii) The following estimates and asymptotics are fulfilled
\[
\lb{eT} |T_1(\l)-T_1^0(\l)|\le 31\vk^2 e^x,\ \ \
|T(\l)-T^0(\l)|\le 1281\vk^2 e^{x+|y|}, \ \ \l\in\ol{\L_4^\pm},
\]
\[
\lb{asm} T_1(\l)=T_1^0(\l)\lt(1+O(z^{-6})\rt),\ \ \
T(\l)=T^0(\l)\lt(1+O(z^{-6})\rt),\ \ \ |\l|\to\iy,\ \ \l\in\cD_1,
\]
\[
\lb{rT} T_1(\l)={e^z\/4}\lt(1+\a_0(\l)+O(z^{-9})\rt),\ \ \
T(\l)={e^z\/2}\lt(2\cos z+e^{\o_1z}\b_{01}(\l)+e^{\o_2z}\b_{02}(\l)+O(z^{-9})\rt)
\]
as $|\l|\to\iy,|y|<\pi,$ and
\[
\lb{lT1}
T_1(\l)={e^{(1+\o_1){z\/2}}\/4}\lt(2\cosh{(1-\o_1)z\/2}+e^{(1-\o_1){z\/2}}\a_0(\l)
+e^{-(1-\o_1){z\/2}}\a_1(\l)+O(z^{-9})\rt),
\]
\[
\lb{lT} T(\l)={e^{(1+\o_1)z}\/2}\lt(1+\b_{01}(\l)+O(z^{-9})\rt)
\]
as $|\l|\to\iy, x-|y|<\pi$.
\end{lemma}

\no {\bf Remark}.
i)  Using the estimates of the functions $\b_{jk},\wt \f_{kk},\wt v_{jk}$ in $\ol{\L_4^\pm}$ we obtain
the estimates of the entire functions $T_1,T$ in $\C$.

\no ii) The conditions $|\l|\to\iy, x-|y|<\pi$ imply $\Re\l\to -\iy$  and $|\Im \l|\to 0$.
Note that
$|e^{(1+\o_1)z}|\to+\iy$ and
$|e^{(1-\o_1)z}|$ is bounded in \er{lT1}, \er{lT}.

\no {\it Proof.}
i) Identity \er{bnk} yields
$$
|b_{jk}(\l)|\le{\|V\|^2\/32|z|^6},\ \ \ |c_{kj}(\l))|\le{\|V\|^2\/16|z|^6},
\ \ \ \l\in\C_\pm,\ \ \ 0\le k\le j\le 3.
$$
Substituting these estimates into \er{al3}-\er{bt} we have estimates of $\a_j,\b_{jk}$
in \er{efj}, \er{vjk}.

We will prove estimate of $\wt\phi_{kk}$ in \er{efj}. Recall that \er{eW1} give
$\F=I_4+\F_1+\F_2+\wt \F$ where $\F_s(\l)=\F_s=\{\f_{s,jk}\}_{j,k=0}^3,s=1,2$ and $\wt \F=\{\wt \f_{jk}\}_{j,k=0}^3.$
Now we will prove that
\[
\lb{ws} \f_{jk}=\d_{jk}+\f_{1,jk}+\f_{2,jk}+\wt \f_{jk},\ \ \ |\wt \f_{jk}|\le 120\vk^3,
\ \ \ \l\in\ol{\L_4^\pm},\ \ \  j,k=0,1,2,3.
\]
We need some simple estimate from the matrix theory. Let
$A=\{a_{ij}\}_{i,j=0}^3$ be a $4\ts 4$-matrix with the usual
matrix norm $|A|$. We prove that
\[
\lb{emn} \max_{0\le i,j\le 3}|a_{ij}|\le 2|A|.
\]
For each vector $x=\{x_k\}_{k=0}^3\in\R^4$ estimates $|x|_\iy\le |x|\le 2|x|_\iy$ hold,
where we denote $|x|=(\sum|x_k|^2)^{1/2}$ and $|x|_\iy=\max|x_k|$.
Then
\[
\lb{pmn} |Ax|_\iy\le |Ax|\le 2|Ax|_\iy\le 2|A|\sum|x_k|.
\]
Let $|a_{pq}|=\max_{0\le i,j\le 3}|a_{ij}|$. We take $x_k=\d_{kq}$ in \er{pmn}.
Then $\sum_k a_{ik}x_k=a_{iq}$ and we obtain $\max_i|a_{iq}|=|a_{pq}|\le 2|A|$,
which yields \er{emn}.
Estimate \er{emn} together with \er{eW2} gives the last estimate \er{ws}.

Substituting \er{He} into \er{W1} we obtain
\[
\lb{W1a}
\f_{1,jk}=-{\o_j\/4z^3}\int_0^1 ds V(s)
\left\{\begin{array}{ccc}
e^{z(1-s)(\o_j-\o_k)}&,&j\ge k\\
e^{-zs(\o_j-\o_k)}&,&j<k\\
\end{array}\right.  ,
\]
and the identity $\int_0^1 V(t)dt=0$ yields $\f_{1,jj}=0,\ \ j=0,1,2,3$.

In order to complete the proof of \er{efj} we have to prove the identity $\f_{2,jj}=\a_j$,
where $\a_j$ are defined by \er{al0}. Recall that
\[
\lb{iW2a}
\F_2=
{\O\/16z^6}\int_0^1du\int_0^1V(u)V(s)\lt(F(0,u,s)-F(1,u,s)\rt)ds,
\]
$$
F(t,u,s)=H(-u)\O H(t-s)-H(t-u)a(u-s),\ \ \ H=\{h_{jk}\}_{j,k=0}^3,
$$
$$
\O=\diag(\o_0,\o_1,\o_2,\o_3),\ \ \ \ \ 
a=\diag(a_{00},a_{01},a_{02},a_{03}).
$$
Then the diagonal entries of the matrix $F=\{f_{jk}\}$ are given by
\[
\lb{fjj} f_{jj}(t,u,s)=\sum_{k=0}^3\o_kh_{jk}(-u)h_{kj}(t-s)-h_{jj}(t-u)a_{0j}(u-s).
\]
In particular,
$$
f_{jj}(0,u,s)=\sum_{k=0}^3\o_kh_{jk}(-u)h_{kj}(-s)-h_{jj}(-u)a_{0j}(u-s).
$$
Recall the identities \er{He}:
\[
\lb{hij} h_{jk}(t)=
\left\{\begin{array}{ccc}
0&,&j\ge k\\
e^{zt(\o_j-\o_k)}&,&j<k\\
\end{array}\right.,\ \ \ t<0,\ \ \
h_{jk}(t)=
\left\{\begin{array}{ccc}
-e^{zt(\o_j-\o_k)}&,&j\ge k\\
0&,&j<k\\
\end{array}\right.,\ \ \ t\ge0.
\]
Then $f_{jj}(0,u,s)=0$ and we have
\[
\lb{f2j} \f_{2,jj}
=-{\o_j\/16z^6}\iint_{[0,1]^2}V(u)V(s)f_{jj}(1,u,s)duds,\ \ \ j=0,1,2,3.
\]
Let $0\le u,s\le 1$.
Substituting \er{hij} into \er{fjj} we obtain
$$
f_{jj}(1,u,s)=-\sum_{k=j+1}^3\o_ke^{z(1-s+u)(\o_k-\o_j)}+a_{0j}(u-s),\ \ \ j=0,1,2,\ \ \
f_{33}(1,u,s)=a_{03}(u-s).
$$
Identity \er{f2j} gives
$$
\f_{2,jj}=-{\o_j\/16z^6}\iint_{[0,1]^2}V(u)V(s)
\lt(a_{0j}(u-s)-\sum_{k=j+1}^3\o_ke^{z(1-s+u)(\o_k-\o_j)}\rt)duds
$$
$$
={\o_j\/16z^6}\lt(-\iint_{[0,1]^2}V(u)V(s)a_{0j}(u-s)duds+
\sum_{k=j+1}^3\o_ke^{z(\o_k-\o_j)}(b_{jk}+b_{kj})
\rt),\ \ \ j=0,1,2,
$$
$$
\f_{2,33}=
-{\o_3\/16z^6}\iint_{[0,1]^2}V(u)V(s)a_{03}(u-s)duds.
$$
Moreover, substituting \er{dkm}, \er{dk} into \er{ynm} we obtain
$$
a_{0j}(t)=e^{-z\o_jt}a_j(t)=
\sum_{k=0}^{j-1}\o_ke^{zt(\o_k-\o_j)},\ \ \ j=1,2,3,\ \ \ a_{00}(t)=0,\ \ \ t<0,
$$
$$
a_{0j}(t)=-\sum_{k=j}^3\o_ke^{zt(\o_k-\o_j)},\ \ \ \ \ j=0,1,2,3,\ \ \ \ t\ge 0.
$$
Then we obtain
$\f_{2,jj}=\a_j$, which yields \er{efj}.

We will prove \er{vjk}. Recall the identity (see \er{T1T},\er{ws}):
$
v_{jk}=\f_{jj}\f_{kk}-\f_{jk}\f_{kj}$ and $\f_{jk}=\d_{jk}+\f_{1,jk}+\f_{2,jk}+\wt \f_{jk}$.
The last identities and $\f_{1,kk}=0,\f_{2,kk}=\a_k$ give
$$
v_{jk}=(1+\f_{1,jj}+\f_{2,jj}+\wt \f_{jj})
(1+\f_{1,kk}+\f_{2,kk}+\wt \f_{kk})
-(\f_{1,jk}+\f_{2,jk}+\wt \f_{jk})(\f_{1,kj}+\f_{2,kj}+\wt \f_{kj}).
$$
$$
=(1+\a_j+\wt \f_{jj})
(1+\a_k+\wt \f_{kk})
-(\f_{1,jk}+\f_{2,jk}+\wt \f_{jk})(\f_{1,kj}+\f_{2,kj}+\wt \f_{kj}),
$$
which yields
\[
\lb{vmn}
v_{jk}=1+v_{2,jk}+\wt v_{jk},\ \ \ v_{2,jk}
=\a_j+\a_k-\f_{1,jk}\f_{1,kj},
\]
where
\[
\lb{wv} \wt v_{jk}=\f_{2,jj}(\f_{2,kk}+\wt \f_{kk})+\wt \f_{jj}\f_{kk}
-\f_{1,jk}(\f_{2,kj}+\wt \f_{kj})-(\f_{2,jk}+\wt \f_{jk})(\f_{1,kj}+\f_{2,kj}+\wt \f_{kj}).
\]
Recall the estimates \er{eW2}:
$|\F(\l)|\le 4, |\F_1(\l)|\le 2\vk, |\F_2(\l)|\le 4\vk^2,
|\wt\F(\l)|\le 60\vk^3.$
Using \er{emn} we have
$$
|\f_{kj}|\le 2|\F(\l)|\le 8,\ \ \
|\f_{1,kj}|\le 4\vk,\ \ \
|\f_{2,kj}|\le 8\vk^2,\ \ \ |\wt \f_{kj}|\le 120\vk^3,\ \ \ 0\le k,j\le 3.
$$
Substituting these estimates into \er{wv} we obtain estimate \er{vjk}.

We will prove $v_{2,jk}=\b_{jk}$. Identity \er{W1a} provides
\[
\lb{w1}
\f_{1,jk}=-{\o_j\/4z^3}\int_0^1\!\!e^{-zt(\o_j-\o_k)}V(t)dt,\ \ \
\f_{1,kj}=-{\o_k\/4z^3}\int_0^1\!\!e^{z(1-t)(\o_k-\o_j)}V(t)dt,
\]
$0\le j<k\le 3.$ Then
\[
\lb{pw1}
\f_{1,jk}\f_{1,kj}={\o_j\o_k\/16z^6}\int_0^1\!\!e^{-zu(\o_j-\o_k)}V(u)du
\int_0^1\!\!e^{z(1-s)(\o_k-\o_j)}V(s)ds
=c_{kj},
\]
where we have used \er{bnk}.
Substituting \er{pw1} into \er{vmn} we obtain
$
v_{2,jk}=\a_{j}+\a_{k}-c_{kj},
$
which yields $v_{2,jk}=\b_{jk}$. Then \er{vmn} gives the first identity \er{vjk}.

\no ii) Let $\l\in\ol{\L_8^\pm}$. Then \er{T1T}-\er{vjk} imply
$$
|T_1(\l)-T_1^0(\l)|\le{1\/4}\sum_0^3|e^{z\o_k}||\a_k(\l)+\wt\phi_{kk}(\l)|
\le e^{x}\max_k|\a_k(\l)+\wt\phi_{kk}(\l)|,
$$
$$
|T(\l)-T^0(\l)|\le{1\/2}\sum_{0\le j<k\le3}|e^{z(\o_j+\o_k)}||\b_{jk}(\l)+\wt v_{jk}(\l)|
\le 3 e^{x+|y|}\max_{0\le j<k\le3}|\b_{jk}(\l)+\wt v_{jk}(\l)|,
$$
which yield \er{eT}. Asymptotics \er{asm} follows from \er{eT}.

We will prove \er{rT}, \er{lT} for $y\ge 0$.
The proof for $y<0$ is similar.
Let $|\l|\to\iy,0\le y<C.$ Then \er{T1T} gives
$$
T_1(\l)={e^z\/4}\lt(\phi_{00}(\l)+O(e^{-x})\rt),\ \ \
T(\l)={e^z\/2}\lt(e^{-iz}v_{01}(\l)+e^{iz}v_{02}(\l)+O(e^{-x})\rt)
$$
and \er{efj}, \er{vjk} yield \er{rT}.
Let $|\l|\to\iy,x-y<\pi.$ Then \er{T1T} gives
$$
T_1={1\/4}\lt(\phi_{00}e^{z}+\phi_{01}e^{-iz}+O(e^{-x})\rt),\ \ \
T(\l)={1\/2}\lt(e^{(1-i)z}v_{01}(\l)+O(1)\rt)
$$
and \er{efj}, \er{vjk} yield \er{lT} and
$$
T_1(\l)={e^{z}\/4}\lt(1+\a_0(\l)+O(z^{-9})\rt)
+{e^{\o_1z}\/4}\lt(1+\a_1(\l)+O(z^{-9})\rt),
$$
which implies \er{lT1}.
$\BBox$

\section {Proof of the main theorems}
\setcounter{equation}{0}

Using the definitions of $\D_1,\D_2, T_1,..$
we obtain following identities
\[
\lb{3f}\lb{ra}
 \D_1^2+\D_2^2=1+T_2, \ \ \ \D_1\D_2=2T_1^2-{T_2+1\/2}={T-1\/2},\ \ \ \r={1-T\/2}+T_1^2,
\]
\[
\lb{5d}
 D_{\pm}=(T_1\mp 1)^2-\r={(2T_1\mp1)^2-T_2\/2}={T\mp4T_1+1\/2},\ \
\ D_+-D_-=-4T_1.
\]
Then by Lemma 2.1, the functions $\D_1+\D_2,\D_1\D_2,
D_\pm,\r$ are entire and are real on the real line.
We need the results about the function $\r$.
Recall $\r^0=({\cosh z-\cos z\/2})^2,z=\l^{1/4}=x+iy.$

\begin{lemma}\lb{ar}
 i) For each $V\in L_0^1(\T)$ the function
$\r={1\/2}(T_2+1)-T_1^2$ satisfies
\[
\lb{3a} |\r(\l)-\r^0(\l)|\le3\vk^2e^{2x+\vk},\ \ \ \l\in\C,
\]
\[
\lb{ero} |\r^0(\l)|>{e^{2x}\/4^2},\ \ \ {\rm if}\ \ \
|\l^{1/4}-(1\pm i)\pi n|\ge{\pi\/2\sqrt 2},\ \ \ n\in\Z,
\]
\[
\lb{aro}  \r(\l)=\r^0(\l)\lt(1+O(\l^{-3/2})\rt),\ \
\   |\l|\to\iy,\  \ \ \l\in\cD_1,
\]
\[
\lb{rr} \r(\l)={e^{2z}\/16}\lt(1+2\a_0(\l)+O(z^{-9})\rt),\ \ \
|\l|\to\iy,\ \ \ |y|<\pi,
\]
\[
\lb{rl} \r(\l)={e^{(1+\o_1)z}\/8}\lt(-1+\cosh(1-\o_1)z
+\a(\l)+O(z^{-9})\rt),\ \ \ |\l|\to\iy,\ \ \ x-|y|<\pi,
\]
where
\[
\lb{an}\a=(1+e^{(1-\o_1)z})\a_0+(1+e^{-(1-\o_1)z})\a_1-2\b_{01},\ \ \
\a(-4(\pi n)^4)={|\hat V_n|^2\/(2\pi n)^6}.
\]

\no ii) Let $V\in L_0^1(\T)$. Then for each integer
$N>\|V\|^{1/3}$ the function $\r(\l)$ has exactly $2N+1$
roots, counted with multiplicity, in the disk
$\{\l:|\l|<4(\pi(N+{1\/2}))^4\}$ and for each $n>N$, exactly two
roots, counted with multiplicity, in the domain
$\{\l:|\l^{1/4}-\pi(1+i)n|<\pi/4\}$. There are no other roots.

\end{lemma}

\no {\it Proof.} i) By Lemma 2.2, $\r$ is entire and real analytic.
The identity $\r={T_2+1\/2}-T_1^2$ yields
$$
|\r-\r^0|=\lt|{T_2-T_2^0\/2}-T_1^2+(T_1^0)^2\rt|
\le{1\/2}\lt|{T_2}-T_2^0\rt|+\lt|T_1-T_1^0\rt|\lt|T_1+T_1^0\rt|.
$$
Then estimates \er{ep}, \er{2.2a} provide
\er{3a}.

Using the identity $\r^0(\l)=-\sinh^2{(1-i)z\/2}\sin^2{(1-i)z\/2},$
and the  estimate $e^{|y|}<4|\sin z|, |z-\pi
n|\ge{\pi\/4},n\in\Z$ (see [PT]), we obtain
$$
|\r^0(\l)|>{1\/16}e^{2|\Im{(1-i)z\/2}|+2|\Im{i(1-i)z\/2}|}
={1\/16}e^{|y+x|+|y-x|} ={e^{2x}\/16},
$$
which yields \er{ero}. Asymptotics \er{aro} follows immediately from \er{3a},\er{ero}.

Substituting \er{rT} into \er{ra} we obtain \er{rr}.
Substituting \er{lT1}, \er{lT} into \er{ra} we obtain
\er{rl}. We prove the second identity \er{an}.
 Then the first identity \er{an} and \er{bnk},\er{bt} give at 
 $\l=-4(\pi n)^4$, i.e. $z=(1+i)\pi n$,
$$
\a(\l)=2(\a_0(\l)+\a_1(\l)-\b_{01}(\l))=2c_{10}(\l)=
{-i\/8z^6}e^{z(-i-1)}( b_{01}(\l)+b_{10}(\l))={|\hat V_n|^2\/(2\pi n)^6}.
$$
ii) Introduce the contour
$C_n(r)=\{\l:|\l^{1/4}-\pi(1+i)n|=\pi r\}$. Let $N_1>N$ be another
integer. Consider the contours
$C_0(N+{1\/2}),C_0(N_1+{1\/2}),C_n({1\/4}),n>N$. Note that we have
$
|\l|\ge 4\pi^4\|V\|_1^{4/3},\ \ \ \vk\le{1\/2\sqrt 2\pi^3}
$
on all contours.
Then  \er{3a}, \er{ero} yield on all contours
$$
|\r(\l)-\r^0(\l)|\le
3\vk^2e^{2x+\vk}\le{e^{2x}\/16}<|\r^0(\l)|.
$$
Hence, by the Rouch\'e theorem, $\r$ has as many roots,
counted with multiplicity, as $\r^0$ in each of
the bounded domains and the remaining unbounded domain. Since
$\r^0(\l)$ has exactly one simple root at $\l=0$ and exactly one root of multiplicity 2 at $-4(\pi n)^4,n\ge 1$,
and since $N_1>N$ can be chosen arbitrarily large,
the point ii) follows.
$\BBox$

Recall that the set $\{\l:D_+(\l)=0\}$ is a periodic spectrum and
the set $\{\l:D_-(\l)=0\}$ is an anti-periodic spectrum. Now we
prove a result about the number of periodic and anti-periodic
eigenvalues in the large disc. Recall
$D_\pm^0=(\cos z\mp 1)(\cosh z\mp 1)$.

\begin{lemma} \lb{Dpm}
Let  $V\in L_0^1(\T)$. Then the following estimates and properties are fulfilled:
\[
\lb{3i}
 |D_\pm(\l)-D_\pm^0(\l)|\le 866\vk^2 e^{x+|y|} ,\ \ \
\l\in\ol{\L_4^\pm}.
\]
\no i) For each integer $N>\|V\|^{1/3}$ the function $D_+$ has
exactly $2N+1$ roots, counted with multiplicity, in the domain
$\{|\l|^{1/4}<2\pi(N+{1\/2})\}$ and for each $n>N$, exactly two
roots, counted with multiplicity, in the domain $\{|\l^{1/4}-2\pi
n|<{\pi\/2}\}$. There are no other roots.

\no ii) For each integer $N>\|V\|^{1/3}$ the function $D_-$ has
exactly $2N$ roots, counted with multiplicity, in the domain
$\{|\l|^{1/4}<2\pi N\}$ and for each $n>N$, exactly two roots,
counted with multiplicity, in the domain $\{|\l^{1/4}-\pi
(2n+1)|<{\pi\/2}\}$. There are no other roots.
\end{lemma}

\no {\it Proof.}
Identities \er{5d} give
$$
|D_\pm-D_\pm^0|\le{|T-T^0|+4|T_1-T_1^0|\/2}
\le\vk^2{1701 e^{x+|y|}+4\cdot 31 e^{x}\/2},
$$
which yields \er{3i}.

\no i) Let $N'>N$ be another integer. Let $\l$ belong to the
contours $C_0(2N+1),C_0(2N'+1),C_{2n}({1\/2}),|n|>N$, where
$C_n(r)=\{\l:|\l^{1/4}-\pi n|=\pi r\},r>0$.  Note that $\vk\le
{1\/(2\pi)^3}$ and
$e^{{1\/2}|y|}<4|\sin{z\/2}|,e^{{1\/2}x}<4|\sinh{z\/2}|,z=\l^{1/4},$
on all contours. Then $e^{{1\/2}(x+|y|)}<16|\sin{z\/2}\sinh{z\/2}|$
and  \er{3i}
 on all contours yield
$$
\lt|D_+(\l)-4\sin^2{z\/2}\sinh^2{z\/2}\rt|\le(30\vk)^2 e^{x+|y|}
<(15\vk)^2 \lt|4\sin{z\/2}\sinh{z\/2}\rt|^2
<\lt|4\sin^2{z\/2}\sinh^2{z\/2}\rt|.
$$
Hence, by Rouch\'e's theorem, $D_+(\l)$ has as many roots, counted
with multiplicities, as $\sin^2{z\/2}\sinh^2{z\/2}$ in each of the
bounded domains and the remaining unbounded domain. Since
$\sin^2{z\/2}\sinh^2{z\/2}$ has exactly one simple root at $\l=0$ and exactly
one root of the multiplicity two
at $(2\pi n)^4,n\ge 1$, and since $N'>N$ can be chosen arbitrarily large, the point i) follows. The proof for $D_-$ is similar.
 $\BBox$

Now we prove our first result about the Lyapunov function $\D=T_1+\sqrt \r$.

\no{\bf Proof of Theorem 1.1.}
By Lemma \ref{ar} ii, the function
$\sqrt{\r}$ is an analytic function in the domain
$\cD_r,r=4\pi^4\|V\|_1^{4/3}$ and it has an analytic continuation onto the two-sheeted Riemann surface. The function $\D$ is analytic on the Riemann surface of the function $\sqrt \r$.
Due to identity \er{2l} all branches of
$\D$ have the forms $\D_m(z)={1\/2}(\t_m(z)+\t_m^{-1}(z)), \ m=1,2$.

\no  i) We prove \er{aD}.
Substituting \er{asm}, \er{aro} into \er{1g} we obtain the
first asymptotics \er{aD}. Substituting \er{asm} and
the first asymptotics \er{aD} into  the identity
$\D_2={T-1\/2\D_1}$ (see \er{3f})
we obtain the second asymptotics in \er{aD}.

By \er{2l}, the matrix $M(\l),\l\in\cD_r,$ for large $r>0$, has the eigenvalues $\t_m(\l)$ satisfying the identities
$\t_m(\l)+\t_m(\l)^{-1}=2\D_m(\l).$
Then $\t_m(\l)$ has the form $
\t_m(\l)=\D_m(\l)+\sqrt{\D_m^2(\l)-1} $, where $\sqrt 1=1$. Asymptotics  \er{aD}
give
$$
\t_1(\l)=\cosh z(1+O(z^{-6}))+\sqrt{\cosh^2 z(1+O(z^{-6}))-1}=e^{z}(1+O(z^{-6})),
$$
$$
\t_2(\l)=\cos z(1+O(z^{-6}))+\sqrt{\cos^2 z(1+O(z^{-6}))-1}=e^{iz}(1+O(z^{-6})),
$$
$ |\l|\to\iy, \l\in\cD_1,$ which yields asymptotics \er{aev}.

\no ii) By the Lyapunov Theorem (see Sect.1), $\l\in \s(\cL)$ iff $\D_m(\l)\in [-1,1]$ for some $m=1,2$. 
 If $\l\in\s(T)$, then $T_1(\l)$ is real. By ii), $\D(\l)$ is
also real. Hence by \er{1g}, $\sqrt{\r(\l)}$ is real and $\r(\l)\ge 0$.

\no iii) Asymptotics from i) yield $\D_1\neq D_2, \t_1\neq t_2$.
Then we have the statements iii) and iv).

v) We have
$
\D_m'={1\/2}(\t_m+\t_m^{-1})'={1\/2}(1-\t_m^{-2})\t_m'\ne 0, m=1,2,
$
since $\t_m\ne 1,\t_m'\ne 0$. 

vi) Let $G_n=(E_n^-,E_n^+)\neq  \es $ for some $n\ge 1$.
It is possible that  $E_n^\pm$ is a periodic or anti-periodic
eigenvalue. Assume that $E_n^+$ is not a  periodic or anti-periodic
eigenvalue. Then $\D_m(E_n^+)\in (-1,1)$ for some $m=1,2$.
If $E_n^+$ is not a branch point, then we have a contradiction.
$\BBox$

We determine the asymptotics of the Lyapunov function
near the positive semi-axis.

\begin{lemma}
Let $V\in L_0^1(\T)$. Then
the following asymptotics are fulfilled
\[
\lb{Dr} \D_1(\l)={e^z\/2}\lt(1+\a_0(\l)+O(z^{-9})\rt),\ \ \
\D_2(\l)=\cos z+{\b(\l)\/2}+O(z^{-9}),\ \ \ 
\]
as $|\l|\to\iy,\ \ \ |y|<\pi$, where
\[
\lb{bn} \b=e^{\o_1z}\b_{01}+e^{\o_2z}\b_{02}-2\a_0\cos z,\ \ \
\b((\pi n)^4)={|\hat V_n|^2\/16(\pi n)^6}.
\]
\end{lemma}

\no {\it Proof}. Substituting \er{rT}, \er{rr} into the identity $\D_1=T_1+\sqrt\r$ we have asymptotics \er{Dr} for $\D_1$. Using identity $\D_2={T-1\/2\D_1}$ and \er{rT} we obtain
$$
\D_2={2\cos z+e^{\o_1z}\b_{01}+e^{\o_2z}\b_{02}+O(z^{-9})\/2(1+\a_0+O(z^{-9}))}
$$
which yields asymptotics \er{Dr} for $\D_2$.
We prove the second identity in \er{bn}. If $\l=(\pi n)^4$, then $z=\pi n$ and we write
$\b=\b((\pi n)^4),b_{jk}=b_{jk}((\pi n)^4),...$ We have
$$
\b=(-1)^n(\b_{02}+\b_{01}-2\a_{0})
=(-1)^n(\a_{1}+\a_{2}-c_{10}-c_{20}).
$$
Identities \er{bnk}-\er{bt} give
$$
\b=(-1)^n(b_{21}+c_{21}+b_{31}+c_{31}-b_{10}+b_{32}+c_{32}-b_{20}-b_{21}-c_{10}-c_{20})
$$
$$
=(-1)^n(b_{21}+b_{12}+b_{31}-b_{10}+b_{32}-b_{20})
+e^{-\pi n}(b_{31}+b_{13}+b_{32}+b_{23}-b_{10}-b_{01}-b_{20}-b_{02}).
$$
We have
$$
b_{21}+b_{12}={|\hat V_n|^2\/16(\pi n)^6},\ \ \
b_{20}+b_{10}=b_{31}+b_{32}
,\ \ \
b_{20}+b_{02}+b_{10}+b_{01}=
b_{31}+b_{13}+b_{32}+b_{23}.
$$
The last identities give the second identity in \er{bn}.$\BBox$

Now we prove the result about the asymptotics of the gaps and the resonance gaps.

\no {\bf Proof of Theorem 1.2.}  i)
Recall that $\{\l_0^+,\l_{2n}^\pm, n\ge 1\}$ is the sequence of zeros of $D_+$ (counted with multiplicity) such that $\l_{0}^+\le \l_{2}^-\le \l_{2}^+\le \l_{4}^-\le\l_{4}^+\le \l_{6}^- \le... $. And $\{\l_{2n-1}^\pm, n\ge 1\}$ is the sequence of zeros of $D_-$ (counted with multiplicity) such that $\l_{1}^-\le \l_{1}^+\le \l_{3}^-\le
\l_{3}^+\le\l_{5}^-\le \l_{5}^+ \le... $.
Lemma \ref{Dpm} gives
that $|(\l_n^{\pm})^{1/4}-\pi n|<{\pi\/2},n>N$ for some $N>0$.
Furthermore, $\l_n^\pm$ are roots of $\D_j^2-1$ for some $j=1,2$.
Asymptotics \er{Dr} of $\D_1$ shows that  $\D_1(\l)>1$ for large $\l>0$. Hence for such $\l$ the spectrum of $\cL$ has multiplicity 2 or 0, and the points $\l_n^\pm$ are roots of $\D_2^2-1$ for $n>N$.

We determine \er{T2-2}. Lemma \ref{Dpm}.ii yields ${\l_{n}^{\pm}}^{1/4}=\pi n+\ve_n, |\ve_n|<{\pi\/2}$  for $n>N$.
Asymptotics \er{Dr} gives
$\D_2(\l_n^\pm)=(-1)^n\cos\ve_n+O(n^{-6}).$
Then the identity $\D_2(\l_n^\pm)=(-1)^n$ gives $\ve_n=O(n^{-3})$.

Now we will improve the asymptotics of $\ve_n$.
Using again \er{Dr} we have
$
(-1)^n\D_2(\l_n^\pm)=\cos \ve_n+{\b(\l_n^\pm)\/2}+O(n^{-9})$.
Note that
$$
\a_{j}(\l_n^\pm)=\a_{j}((\pi n)^4)
+O(1)\max\limits_{|z-\pi n|\le \ve_n}|\a_{j}'(\l)|=\a_{j}((\pi n)^4)+O(n^{-9}),\ \ \ n\to\iy.
$$
since $\a_{j}'(\l)={d\a_{j}\/dz}O(n^{-3})$ and, by \er{bnk},
${d\a_{j}\/dz}=O(n^{-3}), |z-\pi n|\le 1,n\to+\iy$.
The functions $\b_{jk},c_{jk}$ have similar asymptotics.
Then $\b(\l_n^\pm)=\b((\pi n)^4)+O(n^{-9})$ and
$$
\D_2(\l_n^\pm)
=(-1)^n\lt(1-{\ve_n^2\/2}+{\b((\pi n)^4)\/2}+O(n^{-9})\rt),\ \ \
\b((\pi n)^4)={|\hat V_n|^2\/16(\pi n)^6},
$$
see \er{bn}.
Then the identity $\D_2(\l_n^\pm)=(-1)^n$ gives
$ \ve_n^2=\b((\pi n)^4)+O(n^{-9})$ and  we obtain
$$
4(\pi n)^3\ve_n=\pm\sqrt{|\hat V_n|^2+O(n^{-3})}=\pm |\hat V_n|+O(n^{-{3\/2}}), \ \ \ 
$$$$
\l_n^{\pm}=(\pi n+\ve_n)^4=(\pi n)^4+4(\pi n)^3\ve_n+O(n^{-4}),
$$
 which implies \er{T2-2}.

Asymptotics \er{Dr}, \er{T2-2} provide
$-1<\D_2(\l)<1$, as $\l\in(\l_n^+,\l_{n+1}^-)$ and $\D_2(\l)\not\in[-1,1]$, as $\l\in(\l_n^-,\l_n^+),n>N$.
Then each interval $[\l_n^+,\l_{n+1}^-],n>N$ is a spectral band
with multiplicity $2$ and each interval $(\l_n^-,\l_n^+),n>N$ is a gaps.

We will prove \er{T2-3}.
We consider the case $\Im r_n^\pm\ge 0$.
The proof for $\Im r_n^\pm< 0$ is similar.
Lemma \ref{ar}. ii implies
$z=\l^{1/4}=(1+i)\pi n+\d_n,|\d_n|<1$ for $\l=r_n^\pm, n>N$.
Then \er{rl} gives
$
\r(r_n^\pm)={e^{2\pi n}\/8}(\cosh(1+i)\d_n-1+O(n^{-6})).
$
The condition $\r(r_n^\pm)=0$ yields
$\d_n=O(n^{-3})$.

Using \er{rl} again and the asymptotics $\a(r_n^\pm)=\a(-4(\pi n)^4)+O(n^{-9})$
we obtain
$$
\r(r_n^\pm)={e^{2\pi n}\/8}\lt(i\d_n^2+\a(-4(\pi n)^4)+O(n^{-9})\rt).
$$
The condition $\r(r_n^\pm)=0$  yields
$\d_n^2=i\a(-4(\pi n)^4)+O(n^{-9}).$
Identities \er{an} give
$$
\d_n^2={i|\hat V_n|^2+O(n^{-3})\/(2\pi n)^6},\ \ \ \ \ \
\d_n=\pm {(1+i)|\hat V_n|+O(n^{-2/3})\/\sqrt 2(2\pi n)^3},
$$
and then
$$
r_n^{\pm}=\lt((1+i)\pi n+\d_n\rt)^4=-4(\pi n)^4-8(1-i)(\pi n)^3\d_n+O(n^{-4}),
$$
which yields \er{T2-3}.

\no ii) Assume that we have the
periodic spectrum $\l_{0},\l_{2n}^\pm,n\ge 1$.
Using the asymptotics \er{T2-2} and repeating the standard
arguments (see [PT,pp.39-40]) we obtain the Hadamard factorization
$$
D_+(\l)=-{\l-\l_{0}\/4}
\prod_{n\ge 1}{(\l_{2n}^+-\l)(\l_{2n}^--\l)\/(2\pi n)^8}.
$$
By the similar way, we determine $D_-$ by the anti-periodic
spectrum. Using \er{5d} we have $\r$. Thus, we recover the
resonances.

\no iii) Suppose, that we have the periodic spectrum
and the set of the resonances.  Then we determine the functions
$\r$ by the resonances, and $D_+$ by the periodic spectrum. Using
\er{5d} we get $T_1$, $T_2$ and then $D_-$. Thus, we recover the
anti-periodic spectrum. The proof of another case is similar.
$\BBox$

\section {The spectrum for the small potential}
\setcounter{equation}{0}

\no {\bf Proof of Theorem \ref{1.3}.}
Recall that $\cL y=y''''+\g Vy, V\in L_0^1(\T),\g\in\R$ and let $T_m^\g(\l)=T_m(\l,\g V), m=1,2,\r^\g(\l)=\r(\l,\g V),..$.
Due to \er{(4.61)} we have
\[
\lb{Tg} T_m^\g=T_m^0+\g^2T_{m,2}+\g^3\wt\e_m,\ \
|\wt\e_m(\l,\g)|\le {(m\|V\|)^3\/3!|z|^9}e^{m+\vk},
\ \   \l\in\C,\ \
\]
\[
\lb{s(7.1)}
T_{m,2}(\l)={1\/4}\int_0^m
dt\int_0^tV(s)V(t)\vp_{3}^0(m-t+s,\l)\vp_{3}^0(t-s,\l)ds,\ \ \
m=1,2,
\]
where $T_m^0,\r^0$ were given by \er{0}, and $\vp_{3}^0(t,\l)$ was given by \er{6} and $\wt\e_m(\l,\g)$ is a real analytic function
of $(\l,\g)\in\C^2$.  Simple calculations imply
\[
\lb{ze} T_m^0(\l)=1+{m^4\/4!}\l+O(\l^2),\ \ \
\vp_{3}^0(t,\l)={t^3\/6}+O(\l),\ \ \ |\l|\to 0,
\]
uniformly on $t\in[0,2]$.  Substituting this asymptotics into identity \er{s(7.1)} we obtain
\[
\lb{aet}
T_{m,2}(\l)=v_m+O(\l),\ \ \ |\l|\to 0,
\]
where $v_m$ was given by \er{vm}.
Using identity \er{1d} we obtain
\[
\lb{re} \r^\g(\l)
=\r^0(\l)+\g^2\wt\r(\l,\g),\ \ \ 
\r^0(\l)={\l\/4}+O(\l^2),\ \ \ \l\to 0,
\]
\[
\lb{re1} \wt\r(\l,\g)={T_{2,2}\/2}-2T_1^0(\l)T_{1,2}(\l)+O(\g),\ \ \ \g\to 0,
\]
uniformly in any bounded domain in $\C$.
The function $\r^0(\l)$ has simple roots $\l=0$ and
$\r^0(\l),\wt\r(\l,\g)$ are analytic at the points $\l=0,\g=0$.
Applying the Implicit Function Theorem to
$\r^\g(\l)=\r^0(\l)+\g^2 \wt\r(\l,\g)$ and ${\pa\/\pa\l}\r^\g(\l)\rt|_{\l=\g=0}\ne 0$, we
obtain a unique solution $r_{0}^-(\g),|\g|<\ve,r_{0}^-(0)=0$
of the equation $\r^\g(\l)=0,|\g|<\ve$ for some $\ve>0$.

In order to prove asymptotics \er{r0} we rewrite
the equation $\r^\g(\l)=0$ in the form
\[
\lb{erz} {\l\/4}+O(\l^2)=-\g^2 \wt\r(\l,\g),\ \ \ \l=r_0^-(\g),
\]
which yields
$r_0^-(\g)=O(\g^2),\g\to 0$. Then using asymptotics \er{ze}, \er{aet}, \er{re1},
we obtain
\[
\lb{h} \wt\r(\l,\g)={v_2\/2}-2v_1+O(\g+\l).
\]
Substituting the last asymptotics into
\er{erz}, we have \er{r0}.

Identity \er{5d} gives 
$$
 D_+^\g(\l)
=D_+^0(\l)+\g^2 \wt D_+(\l,\g),\ \ \
\wt D_+(\l,\g)=2(2T_1^0-1)T_{1,2}-{T_{2,2}\/2}+O(\g),\ \ \ \g\to 0,
$$
uniformly in any bounded domain in $\C$, and $D_+^0$ was given by \er{D0},
$
D_+^0
=-{\l\/4}+O(\l^2), |\l|\to 0.
$
The function $D_+^0(\l)$ has simple roots $\l=0$ and
$D_+^0(\l),\wt D_+(\l,\g)$ are analytic at the points $\l=0,\g=0$.
Applying the Implicit Function Theorem to
$D_+^\g(\l)=D_+^0(\l)+\g^2 g(\l,\g)$ and ${\pa\/\pa\l}D_+^\g|_{\l=\g=0}\ne 0$, we
obtain a unique solution $\l_{0}^+(\g),|\g|<\ve,\l_{0}^+(0)=0$
of the equation $D_+^\g(\l)=0,|\g|<\ve$ for some $\ve>0$.

We prove asymptotics \er{l0}. We write the equation $D_+^\g(\l)=0$
in the form
\[
\lb{el0} -{\l\/4}+O(\l^2)=-\g^2 \wt D_+(\l,\g),\ \ \ \l=\l_0^+(\g),
\]
which yields $\l_0^+(\g)=O(\g^2),\g\to 0$. Then using \er{ze}, \er{aet}  we obtain
$$
\wt D_+(\l,\g)
=2v_1-{v_2\/2}+O(\g),\ \ \ \g\to 0.
$$
Substituting the last asymptotics into \er{el0}, we have \er{l0}.

We prove \er{bg0}. Substituting asymptotics \er{l0} into \er{Tg} and using \er{ze}, \er{aet} we obtain
\[
\lb{T1g} T_1^\g(\l_0^+)=1-\g^2A+O(\g^3),\ \ \ A={v_2\/12}-{4v_1\/3},\ \ \  \g\to 0.
\]
Using asymptotics \er{r0} we have
$$
\r^\g(\l_0^+)=sy(\g),\ \ y=(\r^\g)'(r_0^-)+O(s),\ \ \
s=\l_0^+-r_0^-\to 0,\ \ as \ \ \g\to 0.
$$
Substituting  $\r^\g(\l_0^+)=sy(\g)$ into the identity $D_+=(T_1-1)^2-\r$ (see \er{5d}),
and using $D_+(\l_0^+)=0$
we obtain
\[
\lb{ers} s=\l_0^+-r_0^-={(T_1^\g(\l_0^+)-1)^2\/y(\g)}.
\]
 Asymptotics \er{re}, \er{re1} and $(\r^0)'(\l)={1\/4}+O(\l),|\l|\to 0$ give
\[
\lb{rg0} y(\g)=(\r^0)'(r_0^-)+O(\g^2)={1\/4}+O(\g^2),
\]
where we have used asymptotics \er{r0}. Substituting \er{T1g}, \er{rg0} into \er{ers} we have \er{bg0}.

Recall the identity  $\D_m^\g=T_1^\g-(-1)^m\sqrt{\r^\g},m=1,2$.  Then
$$
\D_m^\g(\l)=T_1^\g(r_0^-)-(-1)^m \sqrt{\l-r_0^-}\sqrt{y(\g)}+O((\l-r_0^-)^{3\/2}),
\ \ \ \l-r_0^-\to +0.
$$
Hence the function $\D_1^\g$ is increasing in the interval
$(r_0^-,r_0^-+\ve)$ for some $\ve>0$ (see Fig.(\ref{ssp})) and the function $\D_2^\g$ is decreasing in this interval.
Asymptotics \er{T1g} gives
$$
\D_1^\g(r_0^-)=\D_2^\g(r_0^-)=T_1^\g(r_0^-)=1-\g^2A+O(\g^3),\ \ \
 \g\to 0,\ \ \ r_0^-=r_0^-(\g).
$$
Assume that $A>0$.
Then there exists $\d>0$ such that $-1<\D_1^\g(r_0^-)<1$ for each
$\g\in (0,\d)$ and $\D_1^\g$ is increasing in the interval
$(r_0^-,r_0^-+\ve)$ with some $\ve>0$. Then by Theorem 1.1 iv,
$\D_1^\g$ is increasing in the interval $(r_0^-,\l_0(\g))$, where
$\D_1^\g(\l_0(\g))=1$. Hence $\l_0(\g)=\l_{2n}^{\pm}(\g)$ for some
$n$. Note that $\l_0(0)=0$, since $\D_1^0(\l)=\cosh z$. Then
$\l_0(\g)=\l_0^+(\g)=\l_0^+$.
Hence $-1<\D_1^\g(\l)<1,\l\in(r_0^-,\l_0^+)$ and $\D_1^\g(\l_0^+)=1$.
Moreover, substituting identities \er{Tg}, \er{re} into
the identities $\D_2^\g=T_1^\g-\sqrt{\r^\g}$, we obtain
$\D_2^\g=\cos z+o(\g),\g\to 0$. Then the function
$\D_2^\g+1,0\le\g<\d$ has not any zero in the interval $(r_0^-,\l_0^+)$. Then $-1<\D_2^\g(\l)<1$ for each $\l\in(r_0^-,\l_0^+)$. Hence by
Theorem 1.1 ii), the spectral interval $(r_0^-,\l_0^+)$ has
multiplicity 4.

Now we will show that $A>0$ for all $V\in L_0^1(\T),V\ne 0$.
Firstly we prove the following identity
\[
\lb{A} A={1\/2\cdot 288}\int_0^1 f(u)\int_u^1 V(t)V(t-u)dt du,\ \ \
f(u)=u^2(2u^4-6u^3+5u^2-1).
\]
Let $A_0=2\cdot 288 A$.
Identities \er{vm}, \er{bg0} give
\[
\lb{Ah}
A_0={1\/6}\int_0^2 dt\int_0^tV(s)V(t)(2-t+s)^3(t-s)^3ds-{8\/3}h_1
=-{8h_1\/3}+{h_2\/3}+{h_3\/6},
\]
where
$$ \ \ \ \
h_1=\int_0^1 dt\int_0^tV(s)V(t)(1-t+s)^3(t-s)^3ds,
$$
$$
h_2=\int_0^1 dt\int_0^tV(s)V(t)(2-t+s)^3(t-s)^3ds,\ \ \
h_3=\int_0^1 dt\int_0^1V(s)V(t)(1-t+s)^3(1+t-s)^3ds.
$$
Then
\[
\lb{hn} h_1=I_3-3I_4+3I_5-I_6,\ \ \ h_2=8I_3-12I_4+6I_5-I_6,\ \ \ h_3=-6I_2+6I_4-2I_6,
\]
where
\[
\lb{Im} I_m=\int_0^1 dt\int_0^tV(s)V(t)(t-s)^m ds=\int_0^1  u^mdu\int_u^1V(t-u)V(t) dt,\ m\ge 0.
\]
Substituting \er{hn} into \er{Ah} and using \er{Im}  we obtain
\[
\lb{AI} A_0=5I_4-6I_5+2I_6-I_2
=\int_0^1 f(s)\int_s^1 V(t)V(t-s)dt ds,
\]
where
$
f=u^2(2u^4-6u^3+5u^2-1),
$
which yields \er{A}.

Now we prove that $A>0$. Using $f^{(j)}(0)=f^{(j)}(1),0\le j\le 4,
f^{(5)}(1)=-f^{(5)}(0)=6!$ we have
\[
f(t)=\sum_n f_ne^{i2\pi nt},\ \ \ f_n={6!2\/(2\pi n)^6}, \ \ n\ne 0,\ \ \
f_0=-{1\/21},\ \  V(t)=\sum_{n\ne 0}\hat V_ne^{i2\pi nt}.
\]
 Substituting these identities into \er{AI} we get
 $$
A_0=\sum_{m,n} \hat V_n\hat V_m\int_0^1 f(s)ds\int_s^1e^{i2\pi (n+m)t}e^{-i2\pi ns}dt=F_1+F_2
$$
where
$$
F_2=\sum_{m+n\ne 0} {\hat V_n\hat V_m\/2\pi i(n+m)}\int_0^1 f(s)e^{-i2\pi ns}(1-e^{i2\pi (n+m)s})ds
=\sum_{m+n\ne 0} \hat V_n\hat V_m{f_n-f_m\/2\pi i(n+m)}=0
$$
and
$$
F_1=\sum_{-\iy}^\iy |\hat V_n|^2\int_0^1 f(s)(1-s)e^{-i2\pi ns}ds.
$$
Note that
$$
\int_0^1 f(s)(1-s)e^{-i2\pi ns}ds
=\sum_p \int_0^1(1-s) f_pe^{i2\pi (p-n)s}ds=\sum_{p\ne n} {-f_p\/i2\pi (p-n)}+{f_n\/2}.
$$
Then
$$
F_1=\sum_{n\ne 0}|\hat V_n|^2\lt(\sum_{p\ne n}{-f_p\/i2\pi (p-n)}+{f_n\/2}\rt)=\sum_{n\ne 0}|\hat V_n|^2{f_n\/2}>0
$$
since $\hat V_0=0, F_1$ is real, $f_p>0,p\ne 0.\BBox$

\section {Complex resonances.}
\setcounter{equation}{0}

Consider the operator
$\cL^\g={d^4\/dt^4} +\g \d_{per},\g\in\C,$ where $\d_{per},
\ \d_{per}=\sum \d(t-n)$.
Let $T_1^\g=T_1(\cdot,\g \d_{per}),\r^\g=\r(\cdot,\g \d_{per}),...$. Recall that
 $\l=z^4$ and for the  case $V=0$ we have
\[
\lb{uq} T_1^0(\l)=c_+(z),\ \ \r^0(\l)=c_-^2(z),\ \
c_\pm(z)={\cosh z\pm\cos z\/2},\ \ \ s_\pm(z)={\sinh z\pm\sin z\/2}.
\]

\begin{lemma}
For the operator $\cL^\g=d^4/dt^4+\g\d_{per}$ the following identities are fulfilled:
\[
\lb{dp} T_1^\g=T_1^0-\g {s_-\/4z^3},\ \ \
\r^\g=\r^0-\g{c_-s_+\/2z^3}+\g^2{s_-^2\/16z^6}.
\]
\end{lemma}

\no {\it Proof.} The solution $y(t)$ of the equation $y''''+V^\g
y=\l y$ and $y', y''$ are continuous  and
$y'''(n+0)-y'''(n-0)=-\g y(n),\ n\in\Z$.
Then the fundamental solutions $\vp_j(t,\l),j=0,1,2,3,$ have
the form
$$
\vp_j(t)= \vp_j^0(t),\ \ \  0\le t <1,\ \ \
\vp_j(t)=\vp_j^0(t)-\g\vp_3^0(t-1)\vp_j^0(1), \ \ \ 1\le t<2,
$$
$$
\vp_j(t)\!\!=\!\!\vp_j^0(t)-\g\vp_3^0(t-1)\vp_j^0(1)
-\g\vp_3^0(t-2)\lt(\vp_j^0(2)-\g\vp_3^0(1)\vp_j^0(1) \rt),\ 2\le
t<3,
$$
here and below we write $\vp_j(t)=\vp_j(t,\l)$. Then
$$
T_1^\g={1\/4}\sum_0^3\vp_j^{(j)}(1)=
{1\/4}\lt(\sum_0^3(\vp_j^0)^{(j)}(1)-\g\vp_3^0(1)\rt)
=T_1^0-{\g\/4}\vp_3^0(1),
$$
which yields the first identity in \er{dp}, and
$$
T_2^\g={1\/4}\sum_0^3\vp_j^{(j)}(2)
=T_2^0-{\g\/4}\sum_0^3(\vp_3^0)^{(j)}(1)\vp_j^0(1)
-{\g\/4}\lt(\vp_3^0(2)-\g(\vp_3^0(1))^2\rt).
$$
Identities \er{(2.4)} give
$\sum_0^3(\vp_3^0)^{(j)}(1)\vp_j^0(1)=\vp_3^0(2)$ and we obtain
$$
T_2^\g=T_2^0-{\g\/2}\vp_3^0(2)+{\g^2\/4}(\vp_3^0(1))^2.
$$
Then
$$
\r^\g={T_2^\g+1\/2}-(T_1^\g)^2
={1\/2}\lt(T_2^0-{\g\/2}\vp_3^0(2)+{\g^2\/4}(\vp_3^0(1))^2+1\rt)
-\lt(T_1^0-{\g\/4}\vp_3^0(1)\rt)^2
$$
$$
=\r^0-{\g\/4}\vp_3^0(2)+T_1^0{\g\/2}\vp_3^0(1)
+\lt({\g\/4}\vp_3^0(1)\rt)^2.
$$
Using
$$
2T_1^0\vp_3^0(1)-\vp_3^0(2)={1\/2z^3}\lt((\cosh z+\cos z)(\sinh
z-\sin z)- (\sinh 2z-\sin 2z)\rt)=-{2c_-s_+\/z^3},
$$
we obtain the second identity in \er{dp}.
$\BBox$

We shall show the existence of complex resonances.
In this case $\g$ is not small parameter. We rewrite $\r^\g$, given by \er{dp}, in the form
\[
\lb{rgm}
 \r^\g(\l)=\lt(c_-(z)-\g{s_+(z)\/4z^3}\rt)^2- {\g^2\sinh z\sin z\/16z^6}
={s_+^2(z)\/(4z^3)^2}(F_+(z)-\g)(F_-(z)-\g),
\]
\[
\lb{gz}  F_{\pm}(z)={4z^3c_-(z)\/s_+(z)\pm \sqrt{u}}, \ \ \
u=\sinh z\sin z,\ \  z\in E_n=(2\pi n,\pi(2n+1))
\]
where $z=\l^{1/4}$ and $\sqrt 1=1$.
The following properties of  $F=F_+$ are fulfilled:
for each $n\ge 1$ the functions $F$ are analytic on the interval $\e_{n}$ and,
$$
F'(z)\to -\iy\ \ {\rm as}\ \  z\to 2\pi n+0,\ \ \ {\rm and} \ \ F'(z)\to
\iy\ \ {\rm as}\ \ z\to(2n+1)\pi-0.
$$
Hence for each $n\ge 1$ there exist points $z_{n}\in
E_{n}$, such that $F(z_n)=\min_{z\in E_{n}}F(z)$.
Then the Taylor expansion of the function $F$ at the point
$z_n$ is given by
\[
\lb{cgz} F(z)=F(z_n)+{(z-z_{n})^{2}\/2}\lt(F''(z_n)+\wt F(z)\rt),\ \ \
\wt F(z)=O(z-z_{n})\ \ {\rm as}\ \ z-z_{n}\to 0.
\]
Moreover, for each fixed $z\in E_n$, we have 
\[
\lb{aF} F(z)=4z^3(1+O(e^{-z/2}))=4z^3(1+O(e^{-{\pi n}})),\ \ \ 
{\rm as}\ \ n\to \iy.
\]
We prove that there exist the real and non-real branch points of the function $\D^\g(\l)$
for some $\g$. The corresponding behavior of the functions $\r^\g(\l)$, $\D^\g(\l)$ is shown by Fig.\ref{esp}.

Now we prove Proposition \ref{1.4}. Let $\g_{n}=F(z_{n})$.We prove that there exists $N>0$ such that for each $n\ge N$
there exist $\ve_n>0$ and the functions $r_{n}^\pm(\g),\g_n-\ve_n<\g<\g_n+\ve_n$,
such that $r_{n}^\pm(\g)$ are zeros of the function
$\r^\g(\l),r_{n}^\pm(0)=z_n^4$. Moreover, the following asymptotics
are fulfilled:
\[
\lb{rdp} r_{n}^\pm(\g)=z_n^4
\pm 4z_n^3\lt({2\n\/F''(z_n)}\rt)^{1\/2}
+O\lt(\n^{{3\/2}}\rt),\ \ \ F''(z_n)>0\ \ {\rm as} \ \n= \g-\g_n\to 0.
\]
\no {\bf Proof of Proposition \ref{1.4}.}
Differentiation  in identity \er{gz} for $F=F_+$  yields
\[
\lb{F'} F'=gF,\ \ \ g={3\/z}+{s_+\/c_-}-h,\ \ h={a'\/a},\ \ a=s_++\sqrt u.
\]
Asymptotics \er{aF} gives $F(z_n)\ne 0, n\ge N$. Then  using identity $F'(z_n)=0$ we obtain $g(z_n)=0$, hence
\[
\lb{ah} h(z_n)={3\/z_n}+{s_+(z_n)\/c_-(z_n)}=1+{3\/z_n}+O(e^{-2\pi n}),\ \ \ n\to +\iy.
\]
We get
\[
\lb{F''} F''=g'F+gF',\ \ \ F''(z_n)=g'(z_n)F(z_n),
\]
and
\[
\lb{g'} g'=g_0-h',\ \ \ \ g_0(z_n)=-{3\/{z_n}^2}+{c_+(z_n)\/c_-(z_n)}-{s_+^2(z_n)\/c_-^2(z_n)}
=O(n^{-2}),\ \ \ n\to +\iy.
\]
Consider $h'$. Differentiating identities \er{F'} we obtain
\[
\lb{h'} h'={a''\/a}-h^2,\ \ a'=c_++{u'\/2\sqrt u},\ \ \
a''=s_-+{u''\/2\sqrt u}-{1\/4\sqrt u}\rt({u'\/\sqrt u}\rt)^2.
\]
Using identity \er{F'} for $a$ we have
\[
\lb{a''}{a''(z_n)\/a(z_n)}=1+O(e^{-\pi n}),\ \ \ n\to +\iy.
\]
Substituting \er{ah}, \er{a''} into \er{h'} we have $h'(z_n)=-{6\/z_n}+O(n^{-2})$.
Then \er{g'} gives $g'(z_n)={6\/z_n}+O(n^{-2})$. Thus \er{aF}, \er{F''} imply
$$
F''(z_n)=24z_n^2(1+O(n^{-1})),\ \ \ n\to+\iy.
$$
Thus, for each $r>0$ there exists $N>0$ such that  $F''(z_n)\ge r$ for all $n\ge N$.

Let $n\ge N$. Substituting $F(z_n)=\g_n$ into \er{cgz} we have
\[
\lb{TF} F(z)=\g_n+{(z-z_n)^2\/2}\lt(F''(z_n)+\wt F(z)\rt),\ \ \ \ \wt F(z)=O(z-z_n),\ \ \
z-z_n\to 0.
\]
There exists $\d>0$ such that the function $\wt F(z)$ is analytic in
the disk $\{|z-z_n|<\d\}$ and $F''(z_n)+\wt F(z)>0$ for $z-z_n<(-\d,\d)$.
Then using \er{TF} we rewrite the equation $F(z)-\g=0$ in the form
\[
\lb{eFg} \F_+(z,\g)\F_-(z,\g)=0,\ \ \
\F_\pm(z,\g)=\sqrt{\n}\mp {z-z_n\/\sqrt 2}\sqrt{F''(z_n)+\wt F(z)}.
\]
Using $(\F_{\pm})'_z(z_n,\g_n)\ne 0$ and applying the Implicit Function Theorem we obtain that $\F_{\pm}(z,\g)$ has
exactly one simple root $z_{\pm}(\g_n+\n)$ in the disk $\{|\n|<\ve_n\}$ for some $\ve_n>0$ such that
\[
\lb{z+} z_{\pm}(\g)=z_n+{\sqrt{2\n}\/\sqrt{F''(z_n)+\wt F(z_{\pm}(\g))}},
\ \ \ \ \ z_{\pm}(\g_n)=z_n.
\]
Thus, the function $F(z)-\g,|\n|<\ve_n$ has exactly two
zeros $z_\pm(\g)$. Then \er{rgm} gives that the function $\r^\g(\l),|\g-\g_n|<\ve_n$ has exactly two
zeros $r_n^\pm(\g)=z_\pm^4(\g)$. Substituting asymptotics \er{cgz} for $\wt F$ into \er{z+}, we obtain \er{rdp}.
$\BBox$

\no {\bf Acknowledgements.}  
E.K. is supported  by  DFG project BR691/23-1.

\no {\bf References}

\no [BBK] Badanin, A; Br\"uning, J; Korotyaev, E. The Lyapunov
function for Schr\"odinger operator  with periodic  $2\ts 2$
matrix potential, preprint 2005

\no [DS] N.~Dunford and J.~T.~Schwartz, Linear
Operators Part II:  Spectral Theory, Interscience, New York, 1988.

\no  [GT1] J. Garnett, E. Trubowitz: Gaps and bands of one dimensional  periodic Schr\"odinger operators. Comment. Math. Helv. 59, 258-312 (1984)

\no  [GT2] J. Garnett, E. Trubowitz: Gaps and bands of one dimensional
 periodic Schr\"odinger operators II. Comment. Math. Helv. 62, 18-37 (1987).

\no [GL] Gel'fand I., Lidskii, V.: On the structure of the regions
of stability of linear canonical
   systems of differential equations with periodic coefficients. (Russian)
   Uspehi Mat. Nauk (N.S.) 10 (1955), no. 1(63), 3--40.

\no [Ge] Gel'fand, I.:
 Expansion in characteristic functions of an equation with periodic coefficients.
 (Russian) Doklady Akad. Nauk SSSR (N.S.) 73, (1950). 1117--1120.

\no [GO] Galunov, G. V.; Oleinik, V. L. Analysis of the dispersion equation for a negative Dirac "comb".  St. Petersburg Math. J.  4  (1993),  no. 4, 707--720

\no  [Ka] T. Kappeler:
Fibration of the phase space for the Korteveg-de-Vries equation.
Ann. Inst. Fourier (Grenoble), 41, 1, 539-575 (1991).

\no  [KK1] P. Kargaev, E. Korotyaev:
The inverse problem for the Hill operator, direct approach.
Invent. Math.  129, no. 3, 567-593(1997)

\no [Ka] T. Kato. Perturbation theory for linear operators.
Springer-Verlag, Berlin, 1995

\no [K1]  E. Korotyaev. The inverse problem for the Hill operator. I
Internat. Math. Res. Notices, 3(1997), 113--125

\no [K2]  E. Korotyaev. Inverse problem and the trace formula for the Hill operator. II Math. Z. 231(1999), no. 2,    345--368

\no [K3]  E. Korotyaev.  Inverse problem for periodic "weighted" operators,  J. Funct. Anal. 170(2000), no. 1, 188--218

\no [K4] E. Korotyaev.  Marchenko-Ostrovki mapping for periodic
Zakharov-Shabat systems,   J. Differential Equations, 175(2001), no. 2, 244--274

\no [K5] E. Korotyaev.  Inverse Problem and Estimates for Periodic
Zakharov-Shabat systems, J. Reine Angew. Math. 583(2005), 87-115

\no [K6] E. Korotyaev. Characterization of the spectrum of Schr\"odinger operators with periodic distributions. Int. Math. Res. Not.  (2003) no. 37, 2019--2031

\no [Kr] M. Krein. The basic propositions of the theory of
$\lambda$-zones of stability of a canonical system of
linear differential equations with periodic coefficients.
In memory of A. A. Andronov, pp. 413--498.
Izdat. Akad. Nauk SSSR, Moscow, 1955.

\no [Ly] Lyapunov, A.: The general problem of stability of motion, 2
nd ed. Gl. Red. Obschetekh. Lit., Leningrad, Moscow, 1935; reprint
Ann. Math. Studies, no. 17, Princeton Univ. Press, Princeton, N.J.,
1947 

\no [MV] Maksudov, F.; Veliev, O.  Spectral analysis of differential operators with periodic matrix coefficients. (Russian) Differentsial'nye Uravneniya 25 (1989), no. 3, 400--409, 547; translation in Differential Equations 25 (1989), no. 3, 271--277

\no  [MO] V. Marchenko, I. Ostrovski:
A characterization of the spectrum
 of the Hill operator. Math. USSR  Sbornik  26, 493-554 (1975).

\no [Mi1] Misura T. Properties of the spectra of periodic and anti-periodic
boundary value problems generated by Dirac operators. I,II, Theor. Funktsii Funktsional. Anal. i Prilozhen, (Russian), 30 (1978), 90-101; 31 (1979), 102-109

\no [Mi2] Misura T. Finite-zone Dirac operators.  Theor. Funktsii
Funktsional. Anal. i Prilozhen, (Russian), 33 (1980), 107-11.

\no [P1]  Papanicolaou, V. The spectral theory of the vibrating periodic beam. Comm. Math. Phys. 170 (1995), no. 2, 359--373.

\no [P2]  V. G. Papanicolaou, The Periodic Euler-Bernoulli
Equation, Transactions of the American Mathematical Society 355, No. 9(2003), 3727--3759.

\no [P3] V. G. Papanicolaou, An Inverse Spectral Result for the
Periodic Euler-Bernoulli Equation, Indiana University Mathematics
Journal,  53, No. 1 (2004), 223--242

\no [PK1]  V. G. Papanicolaou, D. Kravvaritis, An Inverse
Spectral Problem for the Euler-Bernoulli Equation for the Vibrating Beam,
Inverse Problems, 13(1997), 1083--1092.

\no [PK2] V. G. Papanicolaou, D. Kravvaritis, The Floquet
Theory of the Periodic Euler-Bernoulli Equation, Journal of
Differential Equations, 150(1998), 24--41.

\no [PT]  J. P\"oshel, E.Trubowitz.  Inverse spectral theory.
Pure and Applied Mathematics, 130.
Academic Press, Inc., Boston, MA, 1987. 192 pp.

\no [RS]  M. Reed ; B. Simon. Methods of modern mathematical physics. IV. Analysis of operators. Academic Press, New York-London, 1978

\no [YS]  V. Yakubovich; V. Starzhinskii. Linear
differential equations with periodic coefficients. Vol. 1, 2.
Halsted Press [John Wiley \& Sons] New York-Toronto, 1975.

\end{document}